\documentclass[xypic,amscd,syntonly,amssymb,verbatim,12pt]{amsart}

\usepackage{graphicx}
\usepackage[all]{xy}
\usepackage{amssymb}
\usepackage{amsmath}
\usepackage[mathscr]{euscript}

\usepackage{epsfig}
\theoremstyle{plain}
\newtheorem{Thm}{Theorem}
\newtheorem{Prop}[Thm]{Proposition}

\newtheorem{Lem}[Thm]{Lemma}

 \theoremstyle{definition}
\theoremstyle{remark}

\numberwithin{equation}{section}

\begin{document}
 %\title{Coupling class of actions of reductive groups}
 %\title{Vertex operators between general $B$-branes}
  \title{Cohomological vertex operators}

 \author{ ANDR\'{E}S   VI\~{N}A}
\address{Departamento de F\'{i}sica. Universidad de Oviedo.   Avda Calvo
 Sotelo.     33007 Oviedo. Spain. }
 \email{vina@uniovi.es}
%\thanks{This work has been partially supported by Ministerio de Ciencia y
%Tecnolog\'{\i}a, grant FPA2009-11061}

 \maketitle
\begin{abstract}
%Considering the $B$-branes
 %on a Calabi-Yau manifold $X$ acted by a group $G$
 %as objects in the derived category of coherent sheaves, we give a definition of $G$-equivariant
%branes, which generalizes the concept of equivariant sheaves.
 Given a Calabi-Yau manifold  and considering
 the $B$-branes on it as objects in the derived
 category of coherent sheaves, we   identify the vertex operators for strings between two branes with elements of the cohomology groups     of Ext sheaves.
 We define the correlation functions for these general vertex operators. Strings stretching between two coherent sheaves are studied as homological extensions   of the corresponding branes. In this context, we relate strings between different pairs of branes when there are maps  between these pairs.  We  also interpret some  strings with ghost number $k$ as obstructions for lifts or extensions of    strings with ghost number $k-1$.

\end{abstract}

\section {Introduction} \label{S:intro}

In the nonlinear sigma $B$-model on a Calabi-Yau manifold $X$, the local operators  are defined by
 differential forms
 %on $X$
  of type $(0,q)$ with coefficients in the   vector bundles $\bigwedge^pTX$,
  where $TX$ is the holomorphic tangent bundle of $X$ \cite{Witten}.  The BRST operator in this theory is
  the Dolbeault operator $\bar\partial$; so, the spaces $H^q_{\bar\partial}(X,\,\bigwedge^pTX)$ are the groups
  of the BRST cohomology of the $B$-model.

The coupling of that system to gauge fields is carried out through
Chan-Paton factors \cite{Witten0}. In this way, given two
holomorphic vector bundles $V_1\to X$ and $V_2\to X$, an open
string vertex operator for a string stretching from $V_1$ to $V_2$
is given by a $\bar\partial$-closed $(0,q)$-form  with values in
the vector bundle ${\rm Hom}(V_1,\,V_2)$.    More precisely, the
vertex operators for strings between the $B$-branes defined by
$V_1$ and $V_2$ are the elements of the cohomology groups (see
\cite[Sect. 3.2.2]{Aspin}, \cite[page 207]{Aspin-et})
 \begin{equation}\label{(0q)coh}
  H^{0,q}_{\bar\partial}(X,\,{\rm Hom}(V_1,\,V_2)).
 \end{equation}

 In this note, we will deal with $D$-branes of type $B$   in a compact $n$-dimensional K\"ahler manifold $X$ and with  other objects related to them.
 Such a brane  can be
considered as an  object of $D(X)$,  the bounded derived category of coherent
sheaves on $X$ (see   monograph \cite{Aspin-et}, which includes
an exhaustive  list of specific references). Particular $B$-branes are the holomorphic vector bundles on $X$, and (\ref{(0q)coh}) gives the spaces of local operators for strings stretching between two of these particular branes.

Obviously, (\ref{(0q)coh}) makes non sense for strings between general branes. One of the purposes of this article is to ``extend" that formula to the case
%of the vertex operators
 of strings stretching between two arbitrary branes; that is, to give a general definition of the spaces of vertex operators which coincides with  (\ref{(0q)coh})
 when the  branes are locally free sheaves.

  In   Section \ref{S:Vertex},  we carry out the mentioned ``extension". Firstly, we recall the identification of the space (\ref{(0q)coh}) with 
	%a cohomology group of  a sheaf on $X$. That is, using the Dolbeault resolution
 %  of the sheaf of germens of holomorphic
%sections of ${\rm Hom}(V_1,\,V_2)$, we identify the space
%$H^{0,q}_{\bar\partial}(X,\,{\rm Hom}(V_1,\,V_2))$ with
 the $q$th group of cohomology of $X$ with coefficients
%in a sheaf; more precisely,with coefficients
in the sheaf of holomorphic sections of ${\rm Hom}(V_1,\,V_2)$.
 % The  bounded derived category of coherent sheaves on $X$ will be denoted by $D(X)$. 
On the other hand,  if ${\mathcal F}$ and
${\mathcal G}$ are $B$-branes, i.e.,  objects of the category
$D(X)$,  an open string between ${\mathcal F}$ and ${\mathcal G}$
with ghost number $k$ is an element of the Ext group
$\mathit{Ext}^k({\mathcal F},\,{\mathcal G})$
% , where $i+\,\hbox{ghost number of }\,{\mathcal G}-\,\hbox{ghost number of}\,{\mathcal F}$ can be considered as the ghost number of the corresponding strings
 \cite[Sect. 5.2]{Aspin}. As
 $\mathit{Hom}=\mathit{Ext}^0$,  the strings in ${\rm Hom}(V_1,\,V_2)$ have
 %correspond to the particular case when the
 ghost number  $0$. So, the   vertex operators considered in the space  (\ref{(0q)coh})  are for strings with zero ghost number, and we do not have a definition of local operators for strings between general branes with arbitrary ghost number.

  When one considers   two general banes ${\mathcal F}$ and ${\mathcal G}$   and
  strings with   ghost number $k$, as a direct generalization of formula (\ref{(0q)coh}),
 %not necessarily zero,
 we propose
 $$\bigoplus_q H^q\big(X,\,{\mathcal Ext}^k({\mathcal F},\,{\mathcal G})\big),$$
for  the corresponding space of vertex  operators  (see (\ref{bigoplusk})).
 %This   is the direct generalization of the above result of the result mentioned in the preceding paragraph.
 In this way, the vertex operators are regarded as
elements of sheaf cohomology groups.

%When $X$ is a projective variety,
 These spaces of vertex operators
for strings from a locally free sheaf ${\mathcal E}$ to a coherent sheaf ${\mathcal G}$
   are specially simple:
   % If ${\mathcal E}$ is  a locally free ${\mathcal O}$-module,
 the local operators
corresponding to strings with ghost number $k>0$ are trivial (Proposition  \ref{Pro:vertE}). Moreover,
% if ${\mathcal E}$ is  a locally free sheaf, then  
 the strings with ghost
number $k$ from ${\mathcal E}$ to
% other coherent sheaf 
${\mathcal G}$ can be identified with vertex  operators for strings
with ghost number $0$ stretching from ${\mathcal E}$ to ${\mathcal G}$ (Proposition \ref{P:vertex0}).

 When $X$ is a Calabi-Yau manifold of dimension $n$,
there is a holomorphic volume form  which permits   to identify
the space $\bigwedge^nTX$ with the space of holomorphic $n$-forms.
Thus,  in the sigma $B$-model, the correlation functions
corresponding to local operators can be calculated by integration
on $X$ of corresponding differential forms.

 In our  case, when $X$ is a projective Calabi-Yau  variety,
  using the Serre duality in the derived category $D(X)$,  the holomorphic
volume form and the Yoneda product,
  we will define the correlation functions
  for the local operators introduced above (see (\ref{correlation})).
In Proposition \ref{P:equivalencia}, we prove that this definition   generalizes the one
given in \cite[page 208]{Aspin-et} for operators of strings
between holomorphic vector bundles.

We denote with $\mathfrak{Coh}(X)$ the category of coherent sheaves on $X$. The category $D(X)$, as   derived category of $\mathfrak{Coh}(X)$,
is
%a
triangulated
% category
\cite[Sect 1.5]{Kas-Sch}.
 %, \cite[page 386]{Wei}.
Hence, we can consider distinguished triangles in $D(X)$
\begin{equation}\label{d.t.}
{\mathcal B}\overset{u}{\to}{\mathcal C}\overset{v}{\to}{\mathcal
D}\overset{+1}{\to} {\mathcal B}[1].
 %{\mathcal B}\to{\mathcal C}\to {\mathcal D} \to {\mathcal B}[1].
 \end{equation}
  In physical terms, this triangle can be interpreted as a
possible binding of the branes
  ${\mathcal B}$ and ${\mathcal D}$ to form the brane ${\mathcal C}$ \cite[page 368]{Aspin-et}.
   In general,  a string from ${\mathcal B}$ to a brane ${\mathcal
   G}$ does not admit a lift to a string from
    the binding brane
    ${\mathcal C}$ to ${\mathcal G}$. Dually, not all the strings stretching
   from ${\mathcal F}$ to ${\mathcal D}$ can be extended to
   strings to ${\mathcal C}.$ The  obstructions for
   these lifts and extensions are described in Proposition
   \ref{P:noroughly}.

  In Proposition \ref{ExactSeqMathcalExt}, we relate the vertex
  operators for strings
   %such that some of their % whose end points are
   ending on the brane ${\mathcal C}$ and
   the cohomology of other objects of $D(X)$ determined by the
   triangle.  Under additional hypotheses, the result stated in Proposition \ref{ExactSeqMathcalExt}
   adopts the simpler form which appears in Proposition
   \ref{Bind-Branes}.

Section \ref{S:Local}  concerns  branes defined by coherent
sheaves. As the category $\mathfrak{Coh}(X)$ is abelian, it is possible to define in $\mathfrak{Coh}(X)$
the homological concept of extension of an object by other
\cite{Mitchell}. This fact allows us to study the groups of
strings between two coherent sheaves
 %${\mathcal O}$-modules
 ${\mathcal F}$ and ${\mathcal G}$ in terms of extensions,  without resorting to injective resolutions of ${\mathcal G}$. By interpreting the strings as extensions:

(i) We can show easily the definition of the bifunctors
$\mathit{Ext}^k(\,.\, ,\, .\,)$ on  morphisms between coherent
sheaves. The corresponding group homomorphisms   give   relations
among different string spaces.

(ii) We can describe some strings of ghost number $k+1$ as obstructions on
strings with $k$ ghost number.

To illustrate  item (i), let us consider locally free sheaves
%${\mathcal O}$-modules
 ${\mathcal F}_i$, $i=1,2,3$, with ${\mathcal F}_i$   the sheaf of homolorphic sections of the vector bundle $V_i\to X$.
 Then, by Proposition \ref{P:vertex0}, $\mathit{Ext}^k({\mathcal F}_1,\,{\mathcal F}_2)$ is the
$k$th cohomology group of $X$ with coefficients
  in the vector bundle ${\rm Hom}(V_1,\,V_2)$.
  %, assumed that ${\mathcal F}_i$ is the sheaf of homolorphic sections of the vector bundle $V_i\to X$.
   %(see Proposition \ref{P:vertex0}).
    So, given a gauge transformation $V_2\to V_3$, it  induces an obvious
    homomorphism $\mathit{Ext}^k({\mathcal F}_1,\,{\mathcal F}_2)\to \mathit{Ext}^k({\mathcal F}_1,\,{\mathcal F}_3)$.
   However, when the ${\mathcal F}_i$ are general coherent  sheaves the homomorphism between the corresponding
   Ext groups is not so evident. Nevertheless,
   in terms of extensions, the passage from  $\mathit{Ext}^k({\mathcal F}_1,\,{\mathcal F}_2)$ to $\mathit{Ext}^k({\mathcal F}_1,\,{\mathcal F}_3)$ reduces to the construction of a fibred coproduct.

   With respect to item (ii),  given the subsheaf
${\mathcal B}\subset {\mathcal C}$,
 an object ${\mathcal G}$ of $\mathfrak{Coh}(X)$ and a string ${\mathbf
 R}$ with ghost number $p$ between ${\mathcal B}$ and ${\mathcal
 G}$, we construct an string of ${\mathit Ext}^{p+1}({\mathcal
 C}/{\mathcal B},\,{\mathcal G})$ which is the obstruction for an
 extension  of ${\mathbf R}$ to an string stretching between  ${\mathcal C}$
 and ${\mathcal G}$. Dually, some strings in $\mathit{Ext}^{p+1}({\mathcal F},\,{\mathcal
B})$, when ${\mathcal F}$ is a coherent sheaf,  can be regarded as
obstructions for the lift of elements in
$\mathit{Ext}^{p}({\mathcal F},\,{\mathcal C}/{\mathcal B})$ to
strings from ${\mathcal F}$ to ${\mathcal C}$ (see Proposition
\ref{P: directsumand-p}).

 If $X$ is an algebraic variety and it has a positive line bundle, then each coherent sheaf admits a global syzygy. This fact will permit us
 to   describe   the spaces of vertex operators
 % for strings
 as local extensions. 
  More precisely, given the coherent sheaves ${\mathcal F}$ and ${\mathcal G}$, the vertex operators
  for strings between the branes ${\mathcal F}$ and ${\mathcal G}$ are the {\em local} extensions
  of the sheaf ${\mathcal G}$ by
${\mathcal F}$ (see Theorems \ref{ThmH0Extp} and \ref{ThmHqExtp}).
However, when  ${\mathcal F}$ and ${\mathcal G}$ are $D0$ branes
the local extensions  of ${\mathcal G}$ by ${\mathcal F}$   are,
in fact, global; in other words, the spaces of local operators
are isomorphic with the corresponding spaces of strings (see
Example after Theorem \ref{ThmHqExtp}).

\smallskip
To summarize, we enumerate  some  novel points we deal in this paper, and which have not been considered   hitherto  in the literature. 
\begin{enumerate}
\item The introduction of the spaces of vertex operators for strings with ghost number different from $0.$
 \item The  application of cohomological methods for the interpretation of spaces of strings as obstructions to extensions or liftings.
 \item The definition of the correlation functions for cohomological vertex operators.
 \item The interpretation of the vertex operators as local extensions.
\end{enumerate}

\medskip

{\it Acknowledgement.} I thank E. Sharpe for his useful comments and the referee for his suggestions.

%%%%%%%%%%%%%%%%%%%%%%%%%%%%%%%%%%%%%%%%%%%%%%%%%%%%%%%%%%%%%%%%%%%%%%%%%%%%%%%%%%%%%%%%%%%%%%%%%%%%%%%%%%%%%%%%%%%%%%%%%%%%%%%%%%%%%%
%%%%%%%%%%%%%%%%%%%%%%%%%%%%%%%%%%%%%%%%%%%%%%%%%%%%%%%%%%%%%%%%%%%%%%%%%%%%%%%%%%%%%%%%%%%%%%%%%%%%%%%%%%%%%%%%%%%%%%%%%%%%%%%%%%%%%%%%%

\section {Vertex operators} \label{S:Vertex}

We will denote by ${\mathcal O}$  the sheaf of germs of holomorphic functions on $X$,  and we will put
${\mathcal A}^{p,q}$ for the sheaf of germs of differential forms on $X$
of type $(p,q)$.  By the Poincar\'e lemma relative to the operator
$\bar\partial$, one has the following well-known resolution of the
sheaf  ${\mathcal O}$
\begin{equation}\label{resolutionOmega}
0\to{\mathcal O}\to{\mathcal A}^{0,0}\overset{\bar\partial}
{\to}{\mathcal A}^{0,1}\overset{\bar\partial} {\to}{\mathcal
A}^{0,2}\to\dots
\end{equation}

 Let $V_1$ and $V_2$ be   holomorphic vector bundles on $X$. We denote by $V$   the holomorphic  vector bundle   ${\rm Hom} (V_1,\,V_2)$,  and by
${\mathcal O}(V)$ the sheaf of germs of holomorphic sections of
$V$. We have also the corresponding sheaves of $V$-valued
elements
% $$\Omega^p(V)={\mathcal O}(V)\otimes_{\mathcal O}\Omega^p,\;\;\;\;
$$ {\mathcal A}^{0,q}(V)={\mathcal O}(V)\otimes_{\mathcal O}{\mathcal A}^{0,q}.  $$
As the sheaf ${\mathcal O}(V)$ is locally free,  the tensor
product of ${\mathcal O}(V)$ by the resolution
(\ref{resolutionOmega}) gives the resolution
 \begin{equation}\label{resolutionOmega(V)}
0\to{\mathcal O} (V){\longrightarrow}{\mathcal
A}^{0,0}(V)\overset{1\otimes\bar\partial}
{\longrightarrow}{\mathcal
A}^{0,1}(V)\overset{1\otimes\bar\partial}
{\longrightarrow}{\mathcal A}^{0,2}(V){\longrightarrow}\dots
\end{equation}

We put $A^{p,q}(V):=\Gamma(X,\,{\mathcal A}^{p,q}(V))$ for the
space of sections of the corresponding  sheaf. As
(\ref{resolutionOmega(V)}) is a fine resolution of ${\mathcal O}(V)$,
one has
\begin{equation}\label{Hq(X}
 H^q(X,\,{\mathcal O}(V))=h^q(A^{0,\bullet}(V)).
\end{equation}
 That is,
 % In particular, for $V:={\rm Hom}(V_1,\,V_2)$, where   $V_i$ is a holomorphic vector bundle on $X$, $i=1,2$, we have
  % the sheaf of germs of holomorphic functions on $X$.have,
 %For $p=0$, we have,
 \begin{equation}\label{Hq(X1}
 H^q(X,\,{\mathcal O}({\rm Hom}(V_1,\,V_2)))=h^{q}(A^{0,\bullet}({\rm Hom}(V_1,\,V_2))),
  \end{equation}
   where the right hand side is   the $q$th cohomology  object of the complex $A^{0,\bullet}({\rm Hom}(V_1,\,V_2))$; in other words, the space   of vertex operators (\ref{(0q)coh}).

Thus,  the elements of
%the vector space
 $H^q(X,\,{\mathcal O}({\rm Hom}(V_1,\,V_2)))$  can be considered as local operators, and
%represented by a $\bar\partial$-closed  $(0,q)$-form with coefficients in ${\rm Hom}(V_1,V_2)$, we can take  this space as the set  of vertex
the space of vertex operators for a string between $V_1$ and $V_2$  is the following
direct sum of cohomology groups \cite[page 207]{Aspin-et}
 %, \cite[Sect 3.2.2]{Aspin}
\begin{equation}\label{bigoplus}
\bigoplus_q H^q\big(X,\,{\mathcal O}({\rm Hom}(V_1,\,V_2))\big).
 \end{equation}

We  will write   the space (\ref{bigoplus}) in other equivalent
form, which admits a natural generalization to branes
 which are not   locally free sheaves. For this purpose, we   recall some properties of the functor Ext.
  We set ${\mathcal Hom}(\,.\,,\,.\,) $ for the sheaf functor Hom (see \cite[page 87]{{Kas-Sch}})
$$ {\mathcal Hom}(\,.\,,\,.\,):\mathfrak{Coh}(X)^{\rm op}\times\mathfrak{Coh}(X)\to\mathfrak{Sh},$$
 where $\mathfrak{Sh}$ is the category of  sheaves of ${\mathbb C}$-vector spaces  on $X$. It is easy to check that
  \begin{equation}\label{OHom}
  {\mathcal Hom}\big({\mathcal O}(V_1),\,{\mathcal O} (V_2)  \big)={\mathcal O}({\rm Hom}(V_1,\,V_2)   ).
   \end{equation}

 As we said, the bounded derived category of $\mathfrak{Coh}(X)$ will be denoted by $D(X)$, thus one has the derived functor
 $$R{\mathcal Hom}(\,.\,,\,.\,): D(X)^{\rm op}\times D(X)\to D(\mathfrak{Sh}),$$
  where $D(\mathfrak{Sh})$ is the derived category of $\mathfrak{Sh}$. By definition
  ${\mathcal Ext}^k({\mathcal F},\,{\mathcal G})=  H^k R{\mathcal Hom}({\mathcal F},\,{\mathcal G})).$

 On the other hand, we set
  $\mathit{Hom}$  for the corresponding Hom functor of the category
  $\mathfrak{Coh}(X)$; so, denoting with $\mathfrak{Vec}$    the  category of ${\mathbb C}$-vector spaces,
  %is a functor
 \begin{equation}\label{Hom}
 \mathit{Hom}(\,.\,,\,.\,):\mathfrak{Coh}(X)^{\rm
 op}\times\mathfrak{Coh}(X)\to\mathfrak{Vec}.
  \end{equation}
 Its derived functor
 $$R \mathit{Hom}(\,.\,,\,.\,): D(X)^{\rm op}\times D(X)\to D(\mathfrak{Vec}),$$
 determine
 % is denoted denoted by $\mathit{ Ext}^i$ \cite[page 233]{Hart}.
 the Ext groups  \cite[page 194]{Ge-Ma}
 \begin{equation}\label{Extk}
 \mathit{Ext}^k({\mathcal F},\,{\mathcal G})=H^kR\mathit{Hom}({\mathcal F},\,{\mathcal G}).
  \end{equation}
  This is the space of strings with ghost number $k$ between the branes ${\mathcal F}$ and ${\mathcal G}$.
  Since
  \begin{equation}\label{HomDX}
  H^kR\mathit{Hom}({\mathcal F},\,{\mathcal G})={\rm
  Hom}_{D(X)}({\mathcal F},\,{\mathcal G}[k]),
   \end{equation}
   where ${\mathcal
  G}[k]$ is the complex ${\mathcal G}$ shifted by $k$ to the left, the strings
  between two branes
  can be considered as morphisms of the derived category $D(X)$.

Obviously, the ${\mathcal O}(V_i)$ are objects of $D(X)$. The
equality
$$ {\mathcal Hom}\big({\mathcal O}(V_1),\,{\mathcal O} (V_2)  \big)= {\mathcal Ext}^0({\mathcal O}(V_1),\,{\mathcal O}(V_2)),$$
together with (\ref{OHom}), allows us to write
%then,
%the space of vertex operators  the
(\ref{bigoplus}) as
 %can be written
 \begin{equation}\label{bigoplus1}
 \bigoplus_q H^q\big(X,\,{\mathcal Ext}^0({\mathcal O}(V_1),\,{\mathcal O}(V_2))\big).
 \end{equation}
 To sum up, this is space of local operators for strings in
 $$\mathit{Hom}({\mathcal O}(V_1),\, {\mathcal O}(V_2))=\mathit{Ext}^0({\mathcal O}(V_1),\, {\mathcal O}(V_2)).$$

  %$\mathit{Hom}\big({\mathcal O}(V_1),\,{\mathcal O} (V_2)  \big)$ is the space of strings with ghost number $0$ from  ${\mathcal O}(V_1)$ to  ${\mathcal O}(V_2)$.
   Therefore,  the natural generalization of (\ref{bigoplus}) for the
space of vertex operators for strings with ghost number $k$
stretching between the branes ${\mathcal F}$ and ${\mathcal G}$ is
  \begin{equation}\label{bigoplusk}
 \bigoplus_q H^q\big(X,\,{\mathcal Ext}^k({\mathcal F},\,{\mathcal G})\big).
 \end{equation}

In particular, for a string with ghost number $0$   stretching
from ${\mathcal O}$ to ${\mathcal G}={\mathcal O}(\bigwedge^pTX)$,
one has the space of local operators
$$H^q(X,\,{\mathcal Ext}^0({\mathcal O},\,{\mathcal G})).$$
 As ${\mathcal Hom}({\mathcal O},\, .\,)$ is the identity functor, this space coincides with BRST cohomology group of the sigma $B$-model mentioned in Section \ref{S:intro}.

A consequence of the Grothendieck spectral sequence \cite[page 403]{Wei} \cite[page
207]{Ge-Ma} is the known Local-to-Global Ext spectral sequence,
which allows to determine the Ext groups from the sheaves
${\mathcal Ext}$. That is,  given the objects
${\mathcal F}$ and ${\mathcal G}$  in $D(X)$, the     spectral
sequence 
%of the double complex
\begin{equation}\label{s.s.0}
E_2^{p,q}=H^p(X,\,{\mathcal Ext}^q({\mathcal F},\,{\mathcal G}))
 \end{equation}
abuts to $\mathit{Ext}^k({\mathcal F},\,{\mathcal G})$.
 Thus, the spaces  of strings $\mathit{Ext}^k({\mathcal F},\,{\mathcal G})$ are the
limit of the spectral   sequence determined by the space of vertex
operators for strings between ${\mathcal F}$ and ${\mathcal G}$.

\smallskip

%%%%%%%%%%%%%%%%%%%%%%%%%%%%%%%%%%%%%%%%%%%%%%%%%%%%%%%%%%%%%%%%%%%%%%%%%%%%%%%%%%%%%%%%%%%%%%%%%%%%%%%%%%%%%%%%%%%%%%%%%%%%%%%%%%%%%%%%%%%%%

\noindent
\subsection{Correlation functions.}
If  $X$ is a projective variety of dimension $n$, the Serre
functor $\mathscr{S}$ \cite{B-K} in the category $D(X)$ is
$$\mathscr{S}=(\,.\,)\otimes\omega_X[n]:D(X)\to D(X),$$
where $\omega_X$ is the canonical sheaf of $X$ and $[n]$ denotes
the shifting of the complex by $n$ to the left.
%Thus, one has
%bi-functorial isomorphisms
%$$\xi_{{\mathcal F}{\mathcal G}}:{\rm Hom}_{D(X)}({\mathcal F},\,{\mathcal
%G})\to {\rm Hom}_{D(X)}({\mathcal G},\,\mathscr{S}{\mathcal
%F})^*.$$
Thus, for any ${\mathcal F}$, ${\mathcal H}$ objects of $D(X)$,  one has  a perfect pairing
 $$ {\rm Hom}_{D(X)}({\mathcal F},\,{\mathcal H})
   \otimes {\rm Hom}_{D(X)}({\mathcal
 H},\,{\mathcal F}\otimes\omega_X[n])\to {\mathbb C}.$$

Since ${\mathcal F}\otimes\omega_X[n]=
 %{\mathcal F}[n]\otimes\omega_X =
 ({\mathcal F}\otimes\omega_X)[n],$ the above pairing when  ${\mathcal
 H}={\mathcal F}[n]$ gives
$$  F: Ext^n({\mathcal F},\,{\mathcal F}) \otimes {\rm Hom}_{D(X)}({\mathcal
 F},\,{\mathcal F}\otimes\omega_X)\to{\mathbb C}.$$
 %\end{equation}

%If, furthermore,
When  $X$ is a Calabi-Yau variety, as the canonical bundle  of $X$ is trivial,     there is an
$n$-holomorphic form $\Omega$
 which vanishes nowhere. $\Omega$ is unique up to multiplicative constant, and it can be fixed
 % Fixed such a form, for instance 
  imposing
$$\int_X\Omega\wedge\bar\Omega=1.$$ 

  On the other hand, given an open subset $U$ of $X$, we put
  % by putting
  %$\epsilon(U)={\rm id}\otimes\Omega_{|U}$
  %\begin{equation}\label{epsilonperfect}
  $$\epsilon(U):s\in{\mathcal F}(U)\mapsto s\otimes\Omega_{|U}\in{\mathcal F} (U)\otimes_{{\mathcal O}(U)}\omega_X(U).$$
  % \end{equation}
  Then the maps $\epsilon(U)$ determine a
  %we have a
  morphism  $\epsilon:{\mathcal F}\to{\mathcal F}\otimes\omega_X$, which in turn defines an element $\epsilon$ of
  ${\rm Hom}_{D(X)}({\mathcal F},\,{\mathcal F}\otimes\omega_X).$
  % where $K(X)$ is the homotopy category
 So, we have the map
 % constructed a map
\begin{equation}\label{map-t}
t:{\mathit Ext}^n({\mathcal F},\,{\mathcal F})\to {\mathbb C},
 \end{equation}
 where $t(\sigma)$ is 
 %the image of 
 $F(\sigma\otimes\epsilon)$.
 % by the isomorphism (\ref{isomorp}).

%Let ${\mathcal F}$ and ${\mathcal G}$ be objects of $D(X)$. 
 The element $E_{\infty}^{q,p}$ corresponding to the spectral sequence (\ref{s.s.0})
 %associated to the double complex
 %$$E_2^{q,p}=H^q(X,\,{\mathcal Ext}^p({\mathcal F},\,{\mathcal G}))$$
 is a subquotient of $E_2^{q,p}$. That is, there exists a tower of subspaces
 $$0=B_1^{p,q}\subset B_2^{p,q} \subset\dots \subset Z_2^{p,q}\subset Z_1^{p,q}=E_2^{p,q},$$
 and 
 $$E_{\infty}^{p,q}= Z_{\infty}^{p,q}/B_{\infty}^{p,q};\;\;\;  Z_{\infty}^{p,q}=\cap_r Z_{r}^{p,q},\;\;  B_{\infty}^{p,q}=\cup_r B_{r}^{p,q}.$$
  Given the local operator $a\in H^q(X,\,{\mathcal Ext}^p({\mathcal F},\,{\mathcal G}))$, we define
  $$\alpha=\begin{cases} [a]\in E_{\infty}^{p,q},\;\; \hbox{if}\;\; a \in  Z_{\infty}^{p,q} \\
 0 \in E_{\infty}^{p,q}, \;\; \hbox{otherwise.}
  \end{cases} $$
 As the sequence (\ref{s.s.0}) converges
 %$E_r^{q,p}$
  to $\mathit{Ext}^{q+p}({\mathcal F},\,{\mathcal G})$,  the element $\alpha$ is a string with ghost number $p+q$ stretching between ${\mathcal F}$ and ${\mathcal G}$;  $\alpha\in \mathit{Ext}^{q+p}({\mathcal F},\,{\mathcal G})$.

In particular, given the local operators $a_j\in
H^{q_j}(X,\,{\mathcal Ext}^{p_j}({\mathcal F}_{j-1},\,{\mathcal
F}_j))$ for $j=1,\dots, k$, satisfying $\sum_j (q_j+p_j)=n$ and
${\mathcal F}_k={\mathcal F}_0$. Then the Yoneda composite of the
respective $\alpha_j$
 $$\alpha_1\star\dots\star\alpha_k\in \mathit{Ext}^n({\mathcal F}_0,\,{\mathcal
F}_0).$$

 %$\alpha_j\in \mathit{Ext}^{a_j}({\mathcal
%F}_{j-1},\,{\mathcal F}_{j})$ for $j=1,\dots, k$, satisfying
%$a_1+\dots + a_k=n$ and ${\mathcal F}_k={\mathcal F}_0$. The
%Yoneda composite
%$$\alpha_1\star\dots\star\alpha_k\in \mathit{Ext}^n({\mathcal F}_0,\,{\mathcal
%F}_0).$$
 The $k$-correlation function for the local operators
$a_1,\dots,a_k$ can be defined as the complex number
\begin{equation}\label{correlation}
\langle a_1\dots a_k\rangle:=t(\alpha_1\star\dots\star\alpha_k).
 \end{equation}

%%%%%%%%%%%%%%%%%%%%%%%%%%%%%%%%%%%%%%%%%%%%%%%%%%%%%%%%%%%%%%%%%%%%%%%%%%%%%%%%%%%%%%%%%%%%%%%%%%%%%%%%%%%%%%%%%%%%%%%%%%

\smallskip

\subsection{Distinguished triangles and obstructions.}
Given the distinguished triangle (\ref{d.t.}), as the functors
${\rm Hom}_{D(X)}({\mathcal F},\,.\,)$ and ${\rm Hom}_{D(X)}(\,
.\,,{\mathcal G})$ are cohomological functors, the triangle
(\ref{d.t.}) gives rise to the long exact sequences of groups
  \begin{align}\label{longexact1}
    &\dots\to{\mathit  Ext}^p({\mathcal F},\,{\mathcal B})\to {\mathit  Ext}^p({\mathcal F},\,{\mathcal C})
  \to {\mathit  Ext}^p({\mathcal F},\,{\mathcal D})\to{\mathit  Ext}^{p+1}({\mathcal F},\,{\mathcal B})
 \to \dots \\ \label{longexact}
 &\dots\to{  Ext}^p({\mathcal D},\,{\mathcal G})\to {\mathit
Ext}^p({\mathcal C},\,{\mathcal G})
  \to {\mathit  Ext}^p({\mathcal B},\,{\mathcal G})\to{ \mathit Ext}^{p+1}({\mathcal D},\,{\mathcal G})
 \to \dots
  \end{align}
where ${\mathcal F}$ and ${\mathcal G}$ are objects of $D(X)$.
Thus, we have the following proposition.
 \begin{Prop}\label{P:noroughly} Let us assume   that the brane
 ${\mathcal C}$ can decay into the branes ${\mathcal B}$ and ${\mathcal
 D}$, according to (\ref{d.t.}).
 \begin{enumerate}
 \item If $\tau$ is a string between ${\mathcal
 F}$ and ${\mathcal D}$ with ghost number $p$, then ${\mathit
 Ext}^{p+1}({\mathcal F},\,\tau)$ is the obstruction for a lift of
 $\tau$ to a string stretching from ${\mathcal F}$ to ${\mathcal
 C}$.
 \item If $\rho$ is a string in ${\mathit Ext}^p({\mathcal B},\,{\mathcal
 G})$, then ${\mathit
 Ext}^{p+1}({\rho,\,\mathcal G})$ is the obstruction for an extension of
 $\rho$ to a string between  ${\mathcal C}$ and  ${\mathcal
 G}$.
  \end{enumerate}
 \end{Prop}

Roughly speaking, ${\mathit  Ext}^{p+1}({\mathcal D},\,{\mathcal
G})$ computes the strings from ${\mathcal B}$ to ${\mathcal G}$,
with ghost number $p$, which can not be ``extended" to strings
between ${\mathcal C}$ and ${\mathcal G}$. Dually, ${ \mathit
Ext}^{p+1}({\mathcal F},\,{\mathcal B})$ computes the strings
stretching from ${\mathcal F}$ to ${\mathcal D}$, with ghost
number $p$, which can not be ``lifted" to strings between
${\mathcal F}$ and ${\mathcal C}$.

%Thus, the grading in the space of strings between two branes gives
%some more than a direct sum decomposition. It turns out that
% some elements of ghost number $k$  can be considered, in some extent, as
%obstructions to possible ``extensions" (or ``lifts") of   members with ghost number $k-1$.

Given the distinguished
triangle (\ref{d.t.})
%\begin{equation}\label{d.t.}
%{\mathcal B}\overset{u}{\to}{\mathcal C}\overset{v}{\to}{\mathcal D}\overset{+1}{\to} {\mathcal B}[1]
%\end{equation}
 in $D(X)$ and an object ${\mathcal F}$ of this category, we will show that the spaces of vertex operators  $\{ H^q(X,\,{\mathcal Ext}^k({\mathcal F},{\mathcal C}))\}_q$ (with $k$ fixed) can be included in a long exact sequence with other spaces of local operators determined by the  vertices of the triangle.
%we will see study the relations between the spaceof vertex operators 
 %$$H^q(X,\,{\mathcal Ext}^k({\mathcal F},{\mathcal B})),\;\;
% H^q(X,\,{\mathcal Ext}^k({\mathcal F},{\mathcal C}))\; \hbox{ and}\; \, H^q(X,\,{\mathcal Ext}^k({\mathcal F},{\mathcal
% D})).$$

  To deduce that exact sequence, we will prove that  the functor
  $R{\mathcal Hom}({\mathcal F},\,.\,):D(X)\to D(\mathfrak{Sh})$
  %  between the triangulated categories
   is $t$-exact, i.e., it transforms distinguished triangles into distinguished triangles \cite[page 285]{Ge-Ma}.

Returning to the  triangle (\ref{d.t.}), if ${\mathcal P}$ is a
complex in the category of ${\mathcal O}$-modules quasi-isomorphic
to ${\mathcal F}$ with ${\mathcal P}^i$ a projective object,  then
$${\mathcal Ext}^k({\mathcal F},\,{\mathcal B})=h^k({\mathcal Hom}^{\bullet}({\mathcal P},\,{\mathcal B})),$$
where ${\mathcal Hom}^{\bullet}({\mathcal P},\,{\mathcal
  B})$  is   the complex of ${\mathcal O}$-modules, defined by
  \begin{equation}\label{prod_i}
  {\mathcal Hom}^n({\mathcal P},\,{\mathcal B})=\prod_i{\mathcal Hom}({\mathcal P}^i,\,{\mathcal
  B}^{i+n}).
  \end{equation}

(Dually, we could take an injective resolution of ${\mathcal B}$,
such injective resolutions always exist in the category of
${\mathcal O}$-modules).

 Let ${\mathfrak C}$ denote the category  of  complexes of ${\mathcal O}$-modules. With respect to the triangle (\ref{d.t.}),
 % in the category $D(X)$,
 it is not restrictive to assume that ${\mathcal
C}$ and ${\mathcal D}$ are the mappings cylinder and cone of some
morphism $g$ of ${\mathfrak C}$  with domain
${\mathcal B}$; that is,
 \begin{equation}\label{Assumption}
{\mathcal C}={\rm Cyl}(g),\;\; {\mathcal D}={\rm Con}(g),
 \end{equation}
  and the morphisms $u$ and $v$ are the natural
ones (see Proposition 8 in page 256 of \cite{Ge-Ma}).

In the Appendix, we will prove the following proposition.
\begin{Prop}\label{P:Appendix}
 Let $U$ be an open subset of $X$, then
%$${\mathcal Hom}^{\bullet}({\mathcal P},\,{\rm Con}(g))(U)={\rm
%Con}(\widehat g)_{|\,U}, \;\;\;\; {\mathcal Hom}^{\bullet}({\mathcal
%P},\,{\rm Cyl}(g))(U)={\rm Cyl}(\widehat g)_{|\,U},$$
 $${\mathcal Hom}^{\bullet}({\mathcal P},\,{\rm Con}(g))(U)={\rm
 Con}(\widehat g)(U), \; {\mathcal Hom}^{\bullet}({\mathcal P},\,{\rm Cyl}(g))(U)={\rm Cyl}(\widehat g)(U),$$
$\widehat g$ being ${\mathcal Hom}^{\bullet}({\mathcal P},\,g)$.
\end{Prop}

Applying the functor ${\mathcal Hom}^{\bullet} ({\mathcal P}, \, .
\,)$ to
%the restriction to $U$ of
  the triangle (\ref{d.t.}) and taking into account   (\ref{Assumption}),
  %and taking into account (\ref{CylinderHat}) and (\ref{ConeHat}),
  one obtains the sequence
 \begin{align}
 &{\mathcal Hom}^{\bullet} ({\mathcal P},\,{\mathcal B})\to{\mathcal Hom}^{\bullet} ({\mathcal P},\,  {\rm Cyl}(g))  \to {\mathcal Hom}^{\bullet} ({\mathcal P},\,{\rm Con}(g))\to  \notag\\
 &\to{\mathcal Hom}^{\bullet}({\mathcal P},\,{\mathcal B}[1]).
 \notag
  \end{align}
 By Proposition \ref{P:Appendix},  for any open $U$, we have
  \begin{align}
  &{\mathcal Hom}^{\bullet} ({\mathcal P},\,{\mathcal B})(U)\to  {\rm Cyl}(\widehat g)(U)  \to
   {\rm Con}(\widehat g)(U)\to \notag \\
   &{\mathcal Hom}^{\bullet}({\mathcal P},\,{\mathcal
   B}[1])(U)= {\mathcal Hom}^{\bullet}({\mathcal P},\,{\mathcal
   B})[1](U). \notag
   \end{align}

  Hence, one has the following distinguished triangle in the category $D(X)$
 $${\mathcal Hom}^{\bullet}({\mathcal P},\,{\mathcal B})\to {\rm Cyl}(\widehat g)  \to {\rm Con}(\widehat g) \to{\mathcal
 Hom}^{\bullet}({\mathcal P} ,\,{\mathcal B})[1].$$
  We have proved the following proposition.

  \begin{Prop}\label{t-exactitud}
  The functor $R{\mathcal Hom}({\mathcal F},\,.\,)$ is $t$-exact.
  \end{Prop}

  Similarly, the functor $R{\mathcal Hom}({\,.\,,\,\mathcal F})$ are $t$-exact.

\smallskip

  As the functor $H^0$ is cohomological in any triangulated category \cite[page 283]{Ge-Ma}, from the distinguished
  triangle (\ref{d.t.}) in D(X), by Proposition
  \ref{t-exactitud},
  % \ref{ExactSeqMathcalExt},
   one obtains
  the exact sequence
  \begin{equation}\label{dotsto}
  \dots\to{\mathcal Ext}^k({\mathcal F},\,{\mathcal B})\overset{\alpha}{\to}{\mathcal Ext}^k({\mathcal F},\,{\mathcal C})
\overset{\beta}{\to}{\mathcal Ext}^k({\mathcal F},\,{\mathcal
D})\overset{\gamma}{\to}{\mathcal Ext}^{k+1}({\mathcal
F},\,{\mathcal B})\to\dots
 \end{equation}
 in the category
$\mathfrak{Sh}$ of sheaves of ${\mathbb C}$-vector spaces on $X$.

The following proposition yields the mentioned exact sequence, in which the space of vertex operators $\{H^q(X,\,{\mathcal Ext}^k({\mathcal F},\,{\mathcal C}))\}_q$ are involved.
% with the cohomology of $X$ with respect to other objects of $D(X)$.

\begin{Prop} \label{ExactSeqMathcalExt}
Given the distinguished triangle (\ref{d.t.}), the following cohomology
sequence
%of spaces  vertex operators
 is exact
% \begin{equation}\label{aligndots}
 \begin{align}\label{aligndots}
 \dots &\to H^q(X,\,Ker(\beta))\to H^q(X,\,{\mathcal Ext}^k({\mathcal F},\,{\mathcal C}))  \to  H^q(X,\,Im(\beta))\to   \\
&\to  H^{q+1}(X,\,Ker(\beta))\to\dots \notag.
\end{align}
% \end{equation}
where $\beta: {\mathcal Ext}^k({\mathcal F},\,{\mathcal C})
 {\to}{\mathcal Ext}^k({\mathcal F},\,{\mathcal D})$   is the morphism induced by the arrow $v$ of (\ref{d.t.}).
%${\mathcal Ext}^k({\mathcal F},\,{\mathcal C})
%\overset{\beta}{\to}{\mathcal Ext}^k({\mathcal F},\,{\mathcal D}).$
 \end{Prop}

{\it Proof.} From (\ref{dotsto}), one deduces
 the short exact sequence of sheaves
   $$0\to Ker(\beta)\to{\mathcal Ext}^k({\mathcal F},\,{\mathcal C})\to Im(\beta)\to 0.$$
  The sequence (\ref{aligndots}) is   the corresponding long exact cohomology sequence.
  \qed
 % gives rise to   the corresponding long exact cohomology sequence is (\ref{aligndots}).
%\begin{align} \dots &\to H^q(X,\,Ker(\beta))\to H^q(X,\,{\mathcal Ext}^k({\mathcal F},\,{\mathcal C}))  \to  H^q(X,\,Im(\beta))\to \notag \\
%&\to  H^{q+1}(X,\,Ker(\beta))\to\dots. \notag
%\end{align}

\medskip

%%%%%%%%%%%%%%%%%%%%%%%%%%%%%%%%%%%%%%%%%%%%%%%%%%%%%%%%%%%%%%%%%%%%%%%%%%%%%%%%%%%%%%%%%%%%%%%%%%%%%%%%%%%%%%%%%%%%%%%%%%%
%%%%%%%%%%%%%%%%%%%%%%%%%%%%%%%%%%%%%%%%%%%%%%%%%%%%%%%%%%%%%%%%%%%%%%%%%%%%%%%%%%%%%%%%%%%%%%%%%%

\subsection{Vertex operators for strings between locally free sheaves}\label{Sub:Locally}
 If $\mathcal{E}$ is a
 locally free ${\mathcal O}$-module, then one has the following locally free
 resolution  ${\mathcal E}_{\bullet}$ of ${\mathcal E}$, where ${\mathcal E}_i=0$ for $i\ne 0$
  and ${ \mathcal E}_0={\mathcal E},$
 \begin{equation}\label{ResolutionE}
 \dots\to 0\to0\to {\mathcal E}\to {\mathcal E}\to 0,
  \end{equation}
 which can be used to construct the sheaves ${\mathcal Ext}^k({\mathcal E},\,{\mathcal G})$,
 with ${\mathcal G}$ any coherent ${\mathcal O}$-module (see Proposition 6.5 in page 234 of \cite{Hart}).
 Hence,
  $${\mathcal Ext}^p({\mathcal E},\,{\mathcal G})=h^p({\mathcal
  Hom}({\mathcal E}_{\bullet},\,{\mathcal G}))=0,\;\;\;\hbox{for}\; p> 0.$$
  So, we have the proposition.
   \begin{Prop}\label{Pro:vertE}
    Let ${\mathcal E}$ be a locally free ${\mathcal O}$-module and
   ${\mathcal G}$ an arbitrary coherent  ${\mathcal O}$-module. Then,    the space (\ref{bigoplusk}) of vertex operators for a string of $Ext^k({\mathcal
E},\,{\mathcal G})$  is zero, when $k>0$.
 \end{Prop}

Thus, by the proposition, to get nontrivial local operators for
strings with ghost number $k>0$, it is necessary to consider more
general branes, for instance  coherent not locally free sheaves on
$X$. Among those are the quotient ${\mathcal O}/{\mathcal I}$ of
${\mathcal O}$ by an ideal sheaf.

 Under the hypotheses of the proposition, the Local-to-Global spectral sequence is also
  trivial. Then, one has the following result.

 \begin{Prop}\label{P:vertex0}
  Under the hypotheses of Proposition \ref{Pro:vertE},
  % the  string space
$$\mathit{Ext}^k({\mathcal E},\,{\mathcal G})=
H^k(X,\;{\mathcal Hom}({\mathcal E},\,{\mathcal G})).$$
 %with $k>0$,  is zero
 % for any ${\mathcal O}$-module ${\mathcal G}$.
 \end{Prop}

That is, the strings  with ghost number $k$ starting from a
locally free ${\mathcal O}$-module are vertex operators for
strings with ghost number $0$.

\smallskip

\noindent
\underline{\it Distinguished triangles.} Next,  we will consider
Proposition \ref{ExactSeqMathcalExt}, when in the triangle
(\ref{d.t.}) ${\mathcal B}$,  ${\mathcal C}$ and ${\mathcal D}$ are ${\mathcal
O}$-modules.  If ${\mathcal E}$ is a locally free sheaf, by Proposition \ref{P:vertex0}, 
 the exact sequences (\ref{longexact1}) and (\ref{longexact}) can be written in terms of vertex operators. Thus, one obtains the
  following proposition which relates the spaces of vertex operators
for strings which end at the branes ${\mathcal B}$, ${\mathcal
C}$, ${\mathcal D}$ or begin in these branes.
%  such that ${\mathcal B}$ and ${\mathcal D}$ may bind to form ${\mathcal C}$.
\begin{Prop}\label{Bind-Branes}
Let ${\mathcal B}\to{\mathcal C}\to{\mathcal D}\overset{+1}\to{\mathcal B}[1]$ 
 be a distinguished triangle consisting of
${\mathcal O}$-modules, and ${\mathcal E}$ a locally free ${\mathcal O}$-module. Then the
following sequences of vertex operators are exact
\begin{align}
 \dots & \to H^q(X,\,{\mathcal Hom}({\mathcal E},\,{\mathcal B}))\to H^q(X,{\mathcal Hom}({\mathcal E,\,{\mathcal
 C}}))\to \notag \\
 & \to H^q(X,\, {\mathcal Hom}({\mathcal E},\,{\mathcal D}))
 \to H^{q+1}(X,\, {\mathcal Hom}({\mathcal E},\,{\mathcal B}))
\to\dots \notag
 \end{align}
and
\begin{align}
 \dots & \to H^q(X,\,{\mathcal Hom}({\mathcal D},\,{\mathcal E}))\to H^q(X,{\mathcal Hom}({\mathcal C,\,{\mathcal
 E}}))\to \notag \\
 & \to H^q(X,\, {\mathcal Hom}({\mathcal B},\,{\mathcal E}))
 \to H^{q+1}(X,\, {\mathcal Hom}({\mathcal D},\,{\mathcal E}))
\to\dots \notag
 \end{align}
\end{Prop}

\smallskip
 \noindent
 \underline{\it Correlation functions.}
Next, we will deduce the form adopted by the correlation function
 (\ref{correlation}) for strings between locally free ${\mathcal O}$-modules. For this deduction, it is useful the following lemma.
 
 \begin{Lem} \label{LemmaCorrelation} Let $V$ be a holomorphic vector bundle on $X$ and ${\mathcal F}={\mathcal O}(V)$, then the map $t$ defined in (\ref{map-t})   reduces to
 %$$t:\sigma\in {\mathit Ext}^n({\mathcal O}(V),\,{\mathcal O}(V))\mapsto \int_X{\rm tr}(\sigma)\wedge\Omega\in{\mathbb C},\;\;\:
 $$ t(\sigma)=\int_X{\rm tr}(\sigma)\wedge\Omega, $$ 
  $\sigma\in  H^{0,n}_{\bar\partial} (X,\,{\rm End}(V)).$
 %for $\sigma\in {\mathit Ext}^n({\mathcal O}(V),\,{\mathcal O}(V)).$
 \end{Lem}
 
 {\it Proof.}
 % If ${\mathcal F}={\mathcal O}({\rm End}(V))$
If $V$ is a holomorphic vector bundle on $X$, then the Serre duality gives the pairing
$$H^q(X,\,V)\otimes H^{n-q}(X,\,\omega_X\otimes V^*)\to {\mathbb C},\;\;\; \alpha\otimes\beta\mapsto \int_X\alpha\wedge\beta,$$
where in $\alpha\wedge\beta$ is involved the tautological parings on $V_x\otimes V^*_x$; i.e., the trace on ${\rm End}(V_x)$.
 Thus, by Proposition \ref{P:vertex0} together with (\ref{Hq(X1}) the map $t$
 % defined in (\ref{map-t})
  in this case  takes the form
$$\sigma\in H^{0,n}_{\bar\partial} (X,\,{\rm End}(V))={\mathit Ext}^n({\mathcal O} (V),\,{\mathcal O}(V))\mapsto \int_X{\rm tr}(\sigma)\wedge\Omega.$$

\qed

 Let $V_i$, $i=1,\dots,k$  be  holomorphic vector bundles on $X$, we put  ${\mathcal F}_j={\mathcal O}(V_j)$. By Proposition \ref{P:vertex0} together with (\ref{Hq(X1})
 %$$H^{0,q}_{\bar\}}(X,\,{\rm Hom}(V_i,\,V_j))=\mathit{Ext}^q({\mathcal O}(V_i),\,{\mathcal O}(V_j)).$$
  \begin{equation}\label{AuxYoneda}
  H^{0,q}_{\bar\partial}(X,\,{\rm Hom}(V_i,\,V_j))=\mathit{Ext}^q({\mathcal F}_i,\,{\mathcal F}_j).
   \end{equation}
 The composition pairing between the holomorphic vector bundles
 $${\rm Hom}(V_i,\,V_j)\otimes_{\mathcal O} {\rm Hom}(V_j,\,V_r)\to{\rm Hom}(V_i,\,V_r)$$
   gives rise to the Yoneda pairing   \cite[page 707]{G-H}, which, by (\ref{AuxYoneda}),     adopts the following form
 \begin{align}
 H_{\bar\partial}^{0,p}(X,\,{\rm Hom}(V_i,\,V_j) )\otimes_{\mathcal O}  H_{\bar\partial}^{0,q}(X,\,{\rm Hom}&(V_j,\,V_r))\overset{\star}{\longrightarrow} \notag\\
 & H_{\bar\partial}^{0,p+q}(X,\,{\rm Hom}(V_i,\,V_r)). \notag
  \end{align}
% Where in which
 Being
  %in  %$\star$  are
  involved the cup product and the composition of homomorphisms in the definition of $\star$.

 By Proposition  \ref{Pro:vertE}, to have nonzero local operators
 \begin{equation}\label{ajvertex}
 a_j\in H^{q_j}(X,\,{\mathcal Ext}^{p_j}({\mathcal F}_{j-1},\,{\mathcal
F}_j)),
 \end{equation}
 it necessary that $p_j=0$.
 % If ${\mathcal F}_j={\mathcal O}(V_j)$, with $V_j$ an holomorphic
%vector bundle on $X$, then by Proposition \ref{P:vertex0}
% $$H_{\bar\partial}^{0,q}\big(X,\,{\rm Hom}(V_i,\,V_j)  \big)=H^q\big(X,\,{\mathcal Hom}({\mathcal O}(V_i),\,{\mathcal O}(V_j))\big)=\mathit{Ext}^q({\mathcal F}_i,\,{\mathcal F}_j).$$
  And according to Proposition \ref{P:vertex0},
  % in this particular case,
   the $\alpha_j$'s in (\ref{correlation}) are identical to the $a_j$'s. Thus, given $\alpha_i\in Ext^{q_i}({\mathcal F}_{i-1},\,{\mathcal F}_i)$, if $V_0=V_k$ and $\sum q_i=n$, then
  $$\alpha_1\star\dots\star\alpha_k\in Ext^n({\mathcal F}_0,\,{\mathcal F}_0)=H^{0,n}_{\bar\partial}(X,\,{\rm Hom}(V_0,\,V_0)).$$
  
\begin{Prop}\label{P:equivalencia}
 Given the locally free sheaves $\{{\mathcal F}_i={\mathcal O}(V_i)\}_{i=1,\dots,k}$, with $V_0=V_k$ and the vertex operators
 % $a_j$ introduced in (\re 
 $a_j\in H^{q_j}(X,\,{\mathcal Ext}^{0}({\mathcal F}_{j-1},\,{\mathcal
F}_j)),$ such that $\sum q_j=n$. Then the correlation function
 $$\langle a_1\dots a_k\rangle=\int_X {\rm tr}(\alpha_1\wedge\dots\wedge\alpha_k)\wedge\Omega.$$
\end{Prop}

{\it Proof.} It is a direct consequence of the definition (\ref{correlation}) and Lemma \ref{LemmaCorrelation}.
\qed

 Therefore, the  correlation function defined in (\ref{correlation})   generalizes the one given in \cite[page 208]{Aspin-et} for
 vertex operators which are elements in spaces of the form (\ref{Hq(X1}).

\medskip

\noindent
{\bf Remarks.}
 We summarize the most important points considered in this   section. We have introduced the vertex operators on the manifold $X$ as elements of the cohomology groups of objects in the derived category of sheaves. For the case when $X$ is a projective variety, we have defined the correlation functions for the new operators, and we have proved that they generalize the usual ones for local operators associated to strings between locally free sheaves.

When the brane ${\mathcal C}$ may decay into the branes ${\mathcal B}$ and ${\mathcal D}$, we have proved that strings in  ${\mathit Ext}^{p+1}({\mathcal F},\,{\mathcal B})$ are the obstructions  for the lift of strings from ${\mathcal F}$ to ${\mathcal D}$, with ghost number $p$, to strings between ${\mathcal F}$ to ${\mathcal C}$. Dually,  strings in  ${\mathit Ext}^{p+1}({\mathcal D},\,{\mathcal G})$ are the obstructions  for the extension of strings from ${\mathcal B}$ to ${\mathcal G}$, with ghost number $p$, to strings between ${\mathcal C}$ to ${\mathcal G}$.

 We have also studied the relations between the following sets vertex operators 
$$\{H^q(X,\,{\mathcal Hom}({\mathcal E},\,{\mathcal B}))\}_q,\, \{H^q(X,\,{\mathcal Hom}({\mathcal E},\,{\mathcal C}))\}_q,\,\{H^q(X,\,{\mathcal Hom}({\mathcal E},\,{\mathcal D}))\}_q,$$
 when ${\mathcal E}$ is a locally free module.

%%%%%%%%%%%%%%%%%%%%%%%%%%%%%%%%%%%%%%%%%%%%%%%%%%%%%%%%%%%%%%%%%%%%%%%%%%%%%%%%%%%%%%%%%%%%%%%%%%%%%%%%%%%%%%
%%%%%%%%%%%%%%%%%%%%%%%%%%%%%%%%%%%%%%%%%%%%%%%%%%%%%%%%%%%%%%%%%%%%%%%%%%%%%%%%%%%%%%%%%%%%%%%%%%%%%%%%%%%%

\medskip
\section{Strings between coherent sheaves}\label{S:Local}
%{\sc The spaces $H^0(X,\,{\mathcal Ext}^p({\mathcal F},\,{\mathcal G}))$ and $Ext^1({\mathcal F},\,{\mathcal G})$.}
 In this section,    {\em we only consider branes which are objects of the category} $\mathfrak{Coh}(X)$, that is, coherent
  sheaves on $X$.  The category $\mathfrak{Coh}(X)$ has enough injectives;  hence,  the restriction of the functors $\mathit{Ext}^k$, defined in (\ref{Extk}),
  % functors restricted to the abelian category
 to    $\mathfrak{Coh}(X)$ coincides with the classical derived functors of
 %$Hom_{\mathfrak{Coh}(X)}$
 $\mathit{Hom}$  \cite[Corollary 10.7.5]{Wei}. That is, the spaces of  strings (with a given ghost number) between coherent sheaves can be determined calculating the derived functors of $\mathit{Hom}$
 %$Hom_{\mathfrak{Coh}(X)}$,
   % defined on the abelian category $\mathfrak{Coh}(X)$,
    without passing to the derived category.

   As  $\mathfrak{Coh}(X)$  is an
   abelian category, the Ext groups
  %  $Ext^k({\mathcal F},\,{\mathcal G})$
    can be defined and studied in terms of extensions \cite{Mitchell}. We will adopt this point of view  in this section.  In $\mathfrak{Coh}(X)$   it is possible to construct
    % in $\mathfrak{Coh}(X)$
   %one can  construct in it
    exact sequences, fibred products and coproducts. These tools will permit us to relate string spaces
   between different branes and   analyze   the obstructions to ``extensions" and ``lifts" mentioned
    in the preceding section.

 The functor (\ref{Hom}) is not exact. That is, if ${\mathcal B}$
 is a sub-brane of ${\mathcal C}$,
 there are strings with ghost number $0$ between two   branes ${\mathcal B}$ and ${\mathcal G}$,
 %where ${\mathcal B}$ is a sub-brane of ${\mathcal C}$,
  which do not admit
 an extension to a string from ${\mathcal C}$ to ${\mathcal G}$.
 Dually, there are strings with $0$ ghost number
 from ${\mathcal F}$ to the quotient  ${\mathcal C}/{\mathcal B}$ which do not admit a lift to a string starting from
 ${\mathcal F}$. The ``obstructions" to these processes are described
 by  elements of ${\mathit Ext}^1({\mathcal C}/{\mathcal B},\,{\mathcal G})$
 in the first case, and by elements of $\mathit{ Ext}^1({\mathcal F},\,{\mathcal
 B})$ in the dual case. One can say that the functor $\mathit{ Ext}^1$ is the
 track of the inexactitude of the functor $\mathit{Hom}$. Similarly,
 $\mathit{ Ext}^2$ is consequence of the inexactitude of $\mathit{ Ext}^1$, etc.

More precisely, the $\mathit{Ext}^i$ functors as right derived functors of
the functor $\mathit{Hom}$ form a $\delta$-functor \cite[page 205]{Hart}.
Thus, given the short exact sequence $U$ in the abelian category
$\mathfrak{Coh}(X)$ of coherent ${\mathcal O}$-modules
 \begin{equation}\label{shortexact}
 U:\;\;\;0\to{\mathcal B}\to {\mathcal C}\to{\mathcal D}={\mathcal C}/{\mathcal B}\to 0,
 \end{equation}
  one has the  long exact sequences (\ref{longexact1}) and
  (\ref{longexact}).

%%%%%%%%%%%%%%%%%%%%%%%%%%%%%%%%%%%%%%%%%%%%%%%%%%%%%%%%%%%%%%%%%%%%%%%%%%%%%%%%%%%%%%%%%%%%%%%%%%%%%

\subsection{Extensions of coherent sheaves}\label{S:Extensions}
Here, as we have said, one will consider the spaces of strings between two branes as
groups of extensions of coherent ${\mathcal O}$-modules.

 A length $p$ extension of the coherent ${\mathcal O}$-module
${\mathcal G}$ by the coherent sheaf
 %${\mathcal O}$-module
 ${\mathcal F}$ is an exact sequence of ${\mathcal O}$-modules \cite[page
63]{M-L}
\begin{equation}\label{bf S}
{\mathbf S}:\;\;0\to{\mathcal G}\to {\mathcal H}_{p-1}\to{\mathcal
H}_{p-2}\to\dots \to{\mathcal H}_0\to{\mathcal F}\to 0.
 \end{equation}

  Given  $S$ and $S'$ two extensions of length $1$ of ${\mathcal G}$ by ${\mathcal F}$, we say they are equivalent if there exists a morphism between the exact sequences
  ${ S}$ and ${ S'}$ as the one showed in the following diagram
$$\xymatrix{S:\;\; 0\ar[r] & {\mathcal G}\ar[r]\ar[d]_{1} & {\mathcal H}\ar[r]\ar[d]_{\beta} &
  {\mathcal F} \ar[r]\ar[d]^{1} & 0\\
 S':\;\;  0\ar[r] & {\mathcal G}\ar[r] & {\mathcal H}'\ar[r] &
    {\mathcal F} \ar[r] & 0\,.}$$
  It is easy to check that $\beta$ must be an isomorphism;  thus,
  we have defined an equivalence relation
   in the set of extensions of ${\mathcal G}$ by ${\mathcal F}$. The corresponding
  quotient space can be endowed with the structure of abelian group
  \cite{M-L}. This group is denoted by ${\rm Ext}^1({\mathcal
  F},\,{\mathcal G})$ and is called the group of extensions of ${\mathcal
  G}$ by ${\mathcal F}$. The zero element of this group is defined by the extension
  $$0\to{\mathcal G}\to{\mathcal G}\oplus{\mathcal F}\to{\mathcal
  F}\to 0.$$

  Although  ${\rm Ext}^1({\mathcal
  F},\,{\mathcal G})$  and  ${\mathit Ext}^1({\mathcal
  F},\,{\mathcal G})$ are isomorphic, we will maintain the roman typos for the constructions of the Ext groups we will carry out in this subsection, and by  notational consistency, we will put ${\rm Hom}$ for denoting
$\mathit{Hom}$.

%As $\mathfrak{Coh}(X)$ is a full subcategory of $D(X)$, in this
%subsection we will put ${\rm Hom}$ for denoting
%$Hom_{\mathfrak{Coh}(X)}= Hom_{D(X)}$.

 Given an $1$ extension $S$
of ${\mathcal G}$ by ${\mathcal F}$ and a morphism $\gamma:{\mathcal F}'\to{\mathcal F},$ the pullback construction
permits to define an extension $ S\gamma$ of ${\mathcal G}$ by
${\mathcal F}'$ so that the following diagram is commutative
 $$\xymatrix{S\gamma: & 0  \ar[r] & {\mathcal G}\ar[r]\ar[d]_{1} & {\bullet}\ar[r]\ar[d] &
  {\mathcal F}' \ar[r]\ar[d]^{\gamma} & 0\\
  S: & 0 \ar[r] & {\mathcal G}\ar[r] & {\mathcal H} \ar[r] &
    {\mathcal F} \ar[r] & 0\,,}$$
  where the right hand square is the corresponding  cartesian
    square. We have a group homomorphism
    \begin{equation}\label{StSgamma}
S\in{\rm Ext}^1({\mathcal F},\,{\mathcal G})\longrightarrow
S\gamma\in {\rm Ext}^1({\mathcal F}',\,{\mathcal G}).
\end{equation}

If ${\mathcal F}'$, ${\mathcal F}$ and ${\mathcal H}$ are the sheaves of holomorphic sections of the vector
 bundles   $F'$, $F$ and $H$, respectively, then the sheaf represented by $\bullet$ in the diagram  corresponds to
  the vector bundle whose fibre  at the point $x\in X$ is
$$\{(f',\,h)\in F'_x\times H_x\,|\,\gamma_x(f')=\chi_x(h)\},$$
where $\chi$ denotes the morphism $H\to F$.

\smallskip

    Similarly, given $\alpha:{\mathcal G}\to{\mathcal G}'$,
    %by means of the fibred coproduct one
    one can define the extension
    $\alpha S$  of ${\mathcal G}'$ by ${\mathcal F}$

  $$\xymatrix{S: & 0  \ar[r] & {\mathcal G}\ar[r]\ar[d]_{\alpha} & {\mathcal H}\ar[r]\ar[d] &
  {\mathcal F} \ar[r]\ar[d]^{1} & 0\\
 \alpha S: & 0 \ar[r] & {\mathcal G'}\ar[r] & \bullet \ar[r] &
    {\mathcal F} \ar[r] & 0\,,}$$
    where the left hand square is the fibered  coproduct.
   One has a group homomorphism
\begin{equation}\label{StoalphaS}
S\in{\rm Ext}^1({\mathcal F},\,{\mathcal G})\longrightarrow\alpha
S\in {\rm Ext}^1({\mathcal F},\,{\mathcal G'}).
\end{equation}

    The maps  (\ref{StSgamma})  and  (\ref{StoalphaS})  show the functorial
    character of ${\rm Ext}^1(\,.\,,\,.\,)$.  The above operations are associative, in the sense that the extensions $(\alpha
   S)\gamma$ and
   $\alpha (S\gamma)$ are equivalent. Thus, if the branes are related by the morphisms $\alpha$ and $\gamma$ as above,
    we have the following relation between the corresponding spaces of strings
   $$S\in {\rm Ext}^1({\mathcal F},\,{\mathcal G})\longrightarrow
\alpha S\gamma\in {\rm Ext}^1({\mathcal F}',\,{\mathcal G}').$$

On the other hand, fixed $S\in{\rm Ext}^1({\mathcal F},\,{\mathcal G})$, the preceding constructions can be regarded as maps
\begin{align}\label{Extexactsequence}
&\alpha\in {\rm Hom}({\mathcal G},\,{\mathcal G}')\to \alpha S\in
{\rm
Ext}^1({\mathcal F},\,{\mathcal G'}) \\
 &\gamma\in {\rm Hom}({\mathcal
F}',\,{\mathcal F})\to   S\gamma\in {\rm Ext}^1({\mathcal
F}',\,{\mathcal G}).
\end{align}

Given a sub-brane  ${\mathcal B}$ of a brane of ${\mathcal C}$, one has the exact sequence $U$ defined in (\ref{shortexact}).
 The following proposition asserts that a morphism $\rho:{\mathcal B} \to {\mathcal G}$ can be extended to a morphisms defined on ${\mathcal C}$ if the short exact sequence $\rho U$ splits.
 \begin{Prop}\label{P: directsumand}
 Let     ${\mathcal B}$ be a sub-brane of ${\mathcal C}$. The obstruction to an extension of the string
 $\rho\in{\rm Hom}({\mathcal B},\,{\mathcal G})$ to a string stretching from ${\mathcal C}$ to ${\mathcal G}$ is the string
  %with ghost number $1$
 $\rho U\in{\rm Ext}^1({\mathcal C}/{\mathcal B},\,{\mathcal G})$.
 \end{Prop}
  {\it Proof.}   From the long
exact sequence of Ext goups, one obtains  the following exact
sequence
\begin{equation}\label{Exactrecort}
{\rm Hom}({\mathcal C},\,{\mathcal G})\to {\rm Hom}({\mathcal
B},\,{\mathcal G})\to {\rm Ext}^1( {\mathcal C}/{\mathcal
B},\,{\mathcal G}),
 \end{equation}
  where the first arrow is the restriction homomorphism.    The second map is, by (\ref{Extexactsequence}), the correspondence $\rho\mapsto \rho  U$.

  \qed
\smallskip

\noindent
  {\it Example.} If ${\mathcal C}={\mathcal O}$ and
${\mathcal B}$ is an ideal sheaf ${\mathcal I}$ of ${\mathcal O}$.
Then $Z:={\rm Supp}({\mathcal O}/{\mathcal I})$ is an analytic
subvariety of $X$ and we put ${\mathcal O}_Z$ for the coherent
sheaf ${\mathcal O}/{\mathcal I}$.
%Thus,
%$${\rm Ext}^1({\mathcal O}_Z,\,{\mathcal G})={\rm Hom}({\mathcal
%I},\,{\mathcal G})/{\rm Hom({\mathcal O},\mathcal G}).$$

 Let ${\mathcal G}$ be a coherent ${\mathcal O}$-module, then ${\mathcal Hom}({\mathcal
 O},\,{\mathcal G})={\mathcal G}$. By Proposition \ref{P:vertex0},
 $${\rm Ext}^1({\mathcal O},\,{\mathcal G})=H^1(X,\,{\mathcal G}).$$
 Obviously, the restriction of an element of ${\rm Hom}({\mathcal
 O},\,{\mathcal G})=\Gamma(X,\,{\mathcal G})$ determines a string from ${\mathcal I}$ to
 ${\mathcal G}$ with ghost number $0$.

 If $H^1(X,\,{\mathcal G})=0$, then a part of the long exact sequence  of
 Ext is
 %groups gives
 $${\rm Hom}({\mathcal O},\,{\mathcal G})\to {\rm Hom}({\mathcal I},\,{\mathcal
 G})\to{\rm Ext}^1({\mathcal O}_Z,\,{\mathcal G})\to 0.$$
 Thus,
$${\rm Ext}^1({\mathcal O}_Z,\,{\mathcal G})={\rm Hom}({\mathcal
I},\,{\mathcal G})/{\rm Hom({\mathcal O},\mathcal G}).$$
 Roughly
speaking, in ${\rm Ext}^1({\mathcal O}_Z,\,{\mathcal G})$ are
computed
%(module $\Gamma(X,\,{\mathcal G})$)
 the morphisms from ${\mathcal
I}$ to ${\mathcal G}$ (i.e. strings stretching from ${\mathcal I}$
to ${\mathcal G}$ with $0$ ghost number) that do not admit an
extension to morphisms from ${\mathcal O}$ to ${\mathcal G}$, i.
e. which are not defined by global sections of ${\mathcal G}$.
 %(Sheaf of rings ${\mathcal O}_Z$ can contain nilpotent elements, in
%this case, the pair $(Z, \,{\mathcal O}_Z)$ is a scheme).

\smallskip

The result stated in   Proposition \ref{P: directsumand} has the corresponding dual.
  \begin{Prop}\label{P: directsumandual}
 Let     ${\mathcal B}$ be a sub-brane of ${\mathcal C}$. The obstruction to a lift  of the string
 $\tau\in{\rm Hom}({\mathcal F},\,{\mathcal C}/{\mathcal B})$ to a string stretching from ${\mathcal F}$ to ${\mathcal C}$ is the string
  %with ghost number $1$
 $ U\tau\in{\rm Ext}^1({\mathcal F},\,{\mathcal B})$.
 \end{Prop}

\smallskip

    Given the short exact sequences
$$R:\;0\to{\mathcal G}\to {\mathcal H}_1\to{\mathcal K}\to
0,\;\;\;
R':\;0\to{\mathcal K}\to {\mathcal H}_0\to{\mathcal F}\to 0,$$
 the following  exact sequence is called the Yoneda composite of $R$ and $R'$
 $$R\star R':\;0\to{\mathcal G}\to {\mathcal H}_1\to   {\mathcal H}_0\to{\mathcal F}\to 0.$$

 If $S$ and $S'$ are short exact sequences and $\beta$ a morphism
 such that the the Yoneda composite $(S\beta)\star S'$ is defined,
 then there is a morphism of exact sequences
 \begin{equation}\label{Swiches}
 (S\beta)\star S'\to S\star(\beta S'),
 \end{equation}
 which in general is not an isomorphism.

  The long exact sequence (\ref{bf S}) can be written as a
  composition of short exact sequences
  $${\mathbf S}=S_p\star S_{p-1}\star\dots\star S_1,$$
  by decomposing the maps of (\ref{bf S}) in product of  a monomorphism and an epimorphism.

  One says that the sequences
  $${\mathbf S}=S_p \star\dots\star S_1,\;\;\;{\mathbf S}'=S_p' \star\dots\star S_1',$$
  are equivalent if one can be obtained  from the other through
  switches of the form (\ref{Swiches}). The quotient by this equivalence of
  the set of $p$ fold extensions of ${\mathcal G}$ by ${\mathcal
  F}$ can be endowed with structure of abelian group and is
  denoted ${\rm Ext}^p({\mathcal F},\,{\mathcal G})$.
 As it is known, the  groups  ${\rm Ext}^p({\mathcal
F},\,{\mathcal G})$ are isomorphic to the values taken at
$({\mathcal F},\,{\mathcal G})$  by the corresponding derived
functors $Ext^p$  defined in (\ref{Extk}).
% of the functor ${\rm Hom}(\,.\,,\, .\,)$.
 So,  the short exact sequence (\ref{shortexact}) gives rise to long exact sequences of
 ${\rm Ext}^i$ groups as in (\ref{longexact1})-(\ref{longexact}).  Next, we will construct the morphisms of
 these exact sequences  from operations with extensions.

  Let $\alpha:{\mathcal G}\to{\mathcal G}'$ and $\gamma:{\mathcal F}'\to{\mathcal F}$ be the above
   morphisms of branes.  If ${\mathbf S}$ is the extension (\ref{bf S}), then one defines
  $$\alpha {\mathbf S}:=(\alpha S_p) \star\dots\star S_1\in{\rm Ext}^p({\mathcal F},\,{\mathcal G}'),
  \;\;\; {\mathbf S}\gamma:=S_p \star\dots \star
  (S_1\gamma) \in{\rm Ext}^p({\mathcal F}',\,{\mathcal G}).$$
  Thus, the pair $(\gamma,\,\alpha)$ determines the following
  relation between the  corresponding string spaces
  %Let  $({\mathcal F},\,{\mathcal G})$ and   $({\mathcal F}',\,{\mathcal G}')$ be two pairs of branes. Then the morphisms $\alpha:{\mathcal G}\to{\mathcal G}'$ and $\gamma:{\mathcal F}'\to{\mathcal F}$ determines
    \begin{equation}\label{Bifunctor}
    {\mathbf S}\in{\rm Ext}^p({\mathcal F},\,{\mathcal G})\mapsto \alpha({\mathbf S}\gamma)\in {\rm Ext}^p({\mathcal F}',\,{\mathcal G}').
     \end{equation}
   % The bifunctor ${\rm Ext}^p(\,.\,,\,.\,)$ is a bifunctor.

On the other hand, given  the following extension of ${\mathcal
G}$ by ${\mathcal B}$
 \begin{equation}\label{bf S1}
 {\mathbf R}:\;\; 0\to{\mathcal G}\to {\mathcal
E}_{p-1}\to{\mathcal E}_{p-2}\to\dots \to{\mathcal
E}_0\to{\mathcal B}\to 0,
 \end{equation}
composing  (\ref{bf S1}) with (\ref{shortexact}) one obtaines the
exact sequence
$$  {\mathbf R}\star U:\;\; 0\to{\mathcal G}\to {\mathcal
E}_{p-1} \to\dots \to{\mathcal E}_0\to{\mathcal C}\to {\mathcal
D}\to 0,$$ which determines an element of ${\rm
Ext}^{p+1}({\mathcal D},\,{\mathcal G}).$
 The group homomorphism  induced by the map ${\mathbf R}\mapsto {\mathbf
 R}\star U$ is the connecting homomorphism in the long
exact sequence associated to (\ref{shortexact})
 $$\dots\to{\rm Ext}^p({\mathcal D},\,{\mathcal G})\to {\rm Ext}^p({\mathcal C},\,{\mathcal G})
  \to {\rm Ext}^p({\mathcal B},\,{\mathcal G})\to{\rm Ext}^{p+1}({\mathcal D},\,{\mathcal G})
 \to \dots$$

Analogously, the composition of $U$ with    a $p$ extension  of
${\mathcal D}$ by ${\mathcal F}$  gives the connecting
homomorphisms in the long exact sequence
 $$\dots\to{\rm Ext}^p({\mathcal F},\,{\mathcal B})\to {\rm Ext}^p({\mathcal F},\,{\mathcal C})
  \to {\rm Ext}^p({\mathcal F},\,{\mathcal D})\to{\rm Ext}^{p+1}({\mathcal F},\,{\mathcal B})
 \to \dots$$

 Propositions \ref{P: directsumand} and \ref{P: directsumandual} admit the following generalization.
\begin{Prop}\label{P: directsumand-p}
 Let     ${\mathcal B}$ be a sub-brane of ${\mathcal C}$ and $U$ the exact sequence (\ref{shortexact}).
 \begin{enumerate}
 \item The obstruction to an extension of the string
 ${\mathbf R}\in{\rm Ext}^p({\mathcal B},\,{\mathcal G})$ to a string stretching from ${\mathcal C}$ to ${\mathcal G}$ is the string
  %with ghost number $1$
 $${\mathbf R}\star U\in{\rm Ext}^{p+1}({\mathcal C}/{\mathcal B},\,{\mathcal G}).$$
 \item The obstruction to a lift  of the string
 ${\mathbf T}\in{\rm Ext}^p({\mathcal F},\,{\mathcal C}/{\mathcal B})$ to a string stretching from ${\mathcal F}$ to ${\mathcal C}$ is the string
  %with ghost number $1$
 $$ U\star{\mathbf T}\in{\rm Ext}^{p+1}({\mathcal F},\,{\mathcal B}).$$
 \end{enumerate}
 \end{Prop}

%%%%%%%%%%%%%%%%%%%%%%%%%%%%%%%%%%%%%%%%%%%%%%%%%%%%%%%%%%%%%%%%%%%%%%%%%%%%%%%%%%%%%%%%%%%%%%%%%%%%%%%%%%%%%%%%%%%%%%%%%%%%%%%%%

\subsection{Vertex operators and local extensions}

  From now on in this subsection,  {\em we assume that $X$ is an $n$-dimensional algebraic variety and there exists a positive line bundle on $X$}. So, each coherent sheaf ${\mathcal F}$ on $X$
 admits a resolution
 consisting of locally free sheaves \cite{Fulton, Serre} 
 \begin{equation}\label{Resolution}
 0 \rightarrow {\mathcal
E}_n\overset{\partial}{\rightarrow}
\dots\overset{\partial}{\to}{\mathcal
E}_1\overset{\partial}\rightarrow{\mathcal
E}_0\rightarrow{\mathcal F}\rightarrow 0.
 \end{equation}

Given ${\mathcal F}$ and ${\mathcal G}$ two coherent ${\mathcal
O}$-modules,
%sheaves on $X,$
then the sheaf ${\mathcal Ext}^k({\mathcal F},\,{\mathcal G})$ is
the cohomology
  object $h^k({\mathcal Hom}({\mathcal E}_{\bullet},\,{\mathcal
  G}))$ (by the proposition of   \cite{Hart} in page 234 above mentioned).
  We use these sheaves for determining the space (\ref{bigoplusk})
  of the vertex operators for an open string in the  group
$Ext^k({\mathcal F},\,{\mathcal G})$.
% are elements of the space

 The exact sequence (\ref{Resolution}) gives rise to
the short exact sequence
$$0\to {\mathcal K}_p\to{\mathcal E}_{p-1}\to{\mathcal K}_{p-1}\to
0,$$ where ${\mathcal K}_p$ is the kernel of ${\mathcal
E}_{p-1}\to{\mathcal E}_{p-2}$, or equivalently the cokernel of
${\mathcal E}_{p+1}\to{\mathcal E}_{p}$.

 As ${\mathcal Ext}^p({\mathcal F},\,{\mathcal G})$ is the  $p$-th
cohomology
%object $h^p({\mathcal Hom}({\mathcal E}_{\bullet},\,{\mathcal G}))$
of the complex ${\mathcal Hom}({\mathcal E}_{\bullet},\,{\mathcal
G})$, then ${\mathcal Ext}^p({\mathcal F},\,{\mathcal G})$ is the
sheaf associated to the presheaf
$$U\mapsto h^p\big(\Gamma(U,\, {\mathcal Hom}({\mathcal E}_{\bullet},\,{\mathcal G}))   \big),$$
 where
 $$\Gamma(U,\, {\mathcal Hom}({\mathcal E}_{\bullet},\,{\mathcal G}))={\rm Hom}_{\mathcal O|_U}({\mathcal E}_{\bullet|U},\,{\mathcal G}_{|U}).$$

Given a fine enough open covering ${\mathfrak U}=\{U_{\alpha}\}$ of $X$
  and
   $$f\in H^0(X,\,{\mathcal Ext}^p({\mathcal
F},\,{\mathcal G})),$$
 a global section of  ${\mathcal Ext}^p({\mathcal F},\,{\mathcal G})$,
 they determine an element of
$$h^p(\Gamma(U_{\alpha},\, {\mathcal Hom}({\mathcal
E}_{\bullet},\,{\mathcal G}))),$$
 which in turn is the cohomology
class of a cocycle
\begin{equation}\label{cocyclefalpha}
f_{\alpha}\in{\rm Hom}_{\mathcal
O|U_{\alpha}}({\mathcal E}_{p|U_{\alpha}},\,{\mathcal
G}_{|U_{\alpha}}),
\end{equation}
 i.e. satisfying
$f_{\alpha}\circ\partial_{|U_{\alpha}}=0$. Hence, $f_{\alpha}$
%factorizes through the cokernel of
 admits a unique factorization through the cokernel of
$$\partial_{|U_{\alpha}}:{\mathcal E}_{p+1|U_{\alpha}}\to {\mathcal
E}_{p|U_{\alpha}}.$$
  %Thus, there is a unique morphism
  In this way, $f_{\alpha}$ determines a unique morphism
\begin{equation}\label{Hatf}
\widehat f_{\alpha}: {\mathcal K_p}_{|U_{\alpha}}\to {\mathcal
G}_{|U_{\alpha}}.
 \end{equation}

The pushout, i.e. the fibered coproduct, of the morphisms
${\mathcal K}_{p|U_{\alpha}}\to{\mathcal E}_{p-1|U_{\alpha}}$ and
$\widehat f_{\alpha}$ gives rise to the following commutative diagram
in the category of  ${\mathcal O}_{|U_{\alpha}}$-modules, where
the short sequences are exact
\begin{equation}\label{Diagram1}
\xymatrix{0\ar[r] & {\mathcal K}_{p|U_{\alpha}}\ar[d]_{\widehat f_{\alpha}}\ar[r] &{\mathcal E}_{p-1|U_{\alpha}}\ar[r]\ar[d] & {\mathcal K}_{p-1|U_{\alpha}}\ar[r]\ar[d]^{1} & 0 \\
 0\ar[r] & {\mathcal G}_{|U_{\alpha}}\ar[r] & {\mathcal H}^{\alpha}_{p-1}\ar[r] &  {\mathcal K}_{p-1|U_{\alpha}} \ar[r] & 0\; .
}
\end{equation}
 The sequence at the bottom  in this diagram can be jointed with
 the exact sequence
 $$0 \rightarrow {\mathcal
K}_{p-1}\rightarrow{\mathcal E}_{p-2} \rightarrow{\mathcal
E}_{p-3}\rightarrow\dots\rightarrow{\mathcal
E}_0\rightarrow{\mathcal F}\rightarrow 0 $$
 restricted to $U_{\alpha}$ and we obtain a    $p$-extension of ${\mathcal G}_{|U_{\alpha}}$ by ${\mathcal
F}_{|U_{\alpha}}$; i. e. the following exact sequence
 \begin{equation}\label{ExtensionHalpha}
 0 \to{\mathcal G}_{|U_{\alpha}}\to{\mathcal H}^{\alpha}_{p-1}\rightarrow
  {\mathcal E}_{p-2|U_{\alpha}}
\rightarrow\dots\rightarrow{\mathcal
E}_{0|U_{\alpha}}\rightarrow{\mathcal F}_{|U_{\alpha}}\rightarrow
0 .
 \end{equation}

If $U_{\alpha\beta}:=U_{\alpha}\cap U_{\beta}\ne\emptyset$, the
obvious restrictions give rise to the following commutative
diagram
$$\xymatrix{ 0\ar[r] & {\mathcal G}_{|U_{\alpha}}\ar[r]\ar[d] & {\mathcal H}^{\alpha}_{p-1}\ar[r]\ar[d] &\dots\ar[r] &  {\mathcal F}_{|U_{\alpha}} \ar[r]\ar[d] & 0\\
  0\ar[r] & {\mathcal G}_{|U_{\alpha\beta}}\ar[r] & {\mathcal H}^{\alpha}_{p-1|U_{\alpha\beta}}\ar[r] & \dots\ar[r] &   {\mathcal F}_{|U_{\alpha\beta}} \ar[r] & 0\;.
}
$$

 Since the restrictions of $\widehat f_{\alpha}$ and $\widehat f_{\beta}$ to
 $U_{\alpha\beta}$ patch together, one has isomorphisms
 $$\xi_{\beta\alpha}:{\mathcal H}^{\alpha}_{p-1|U_{\alpha\beta}}\to {\mathcal
 H}^{\beta}_{p-1|U_{\alpha\beta}}.$$
 In general,
 $\xi_{\alpha\beta}\circ\xi_{\beta\gamma}\ne\xi_{\alpha\gamma}$.
 Thus, the local extensions (\ref{ExtensionHalpha})
 % ${\mathcal H}^{\alpha}$
  do not define a
 global extension of ${\mathcal G}$ by ${\mathcal F}$. We can
 state the following theorem.

\begin{Thm}\label{ThmH0Extp} Let ${\mathfrak U}=\{U_{\alpha}\}$ be a sufficiently fine open covering
of $X$, ${\mathcal F}$ and ${\mathcal G}$ coherent ${\mathcal
O}$-modules. Each local operator of  $H^0(X,\,{\mathcal
Ext}^p({\mathcal F},\,{\mathcal G}))$, with $p>0$ determines a $p$
extension of ${\mathcal G}_{|U_{\alpha}}$ by ${\mathcal
F}_{|U_{\alpha}}$, for all $\alpha$.
\end{Thm}

\noindent
  {\it Remark 1.} The {\em local} character of spaces  $H^0(X,\,{\mathcal Ext}^p({\mathcal
F},\,{\mathcal G}))$   appears  clearer when they are compared
with the $Ext$ groups. For example, an element of $Ext^1({\mathcal
F},\,{\mathcal G})$ is given by a cocycle of the hypercohomology
${\mathbb H}^1({\mathfrak U},\,{\mathcal Hom}({\mathcal
E}_{\bullet},\,{\mathcal G}))$. This cocycle is a pair
$(f=\{f_{\alpha}\},\,h=\{h_{\alpha\beta}\})$ with
 $$f_{\alpha}\in\Gamma(U_{\alpha},\,{\mathcal Hom}({\mathcal
E}_1,\,{\mathcal
G})),\;\;h_{\alpha\beta}\in\Gamma(U_{\alpha\beta},\,{\mathcal
Hom}({\mathcal E}_0,\,{\mathcal G}))$$ satisfying
 \begin{equation}\label{globaxten}
\partial^*f=0,\;\;\; \delta f=\partial^* h,\;\;\;\delta h=0,
 \end{equation}
 where
$\delta$ denotes the \v{C}ech coboundary operator and $\partial^*$
the one induced by $\partial$ on ${\mathcal Hom}({\mathcal
E}_{\bullet},\,{\mathcal G})$.

The condition $\partial ^*f=0$ implies
$f_\alpha\circ\partial_{|U_{\alpha}}=0$. Hence $f_{\alpha}$
factorizes through the cokernel of $\partial_{|U_{\alpha}}$. Thus,
it is possible to carry out the preceding pushout construction and
we obtain the diagram (\ref{Diagram1}), with $p=1$.  The condition
$\delta f=\partial ^*$ implies that the local extensions
${\mathcal H}^{\alpha}_0$  and ${\mathcal H}^{\beta}_0$ patch together on
%agree when are restricted to
 $U_{\alpha\beta}$. Now,  from the condition
$\delta h=0$, ones
 deduces
% that the isomorphisms $\{\xi_{\alpha\beta}\}$ satisfy
the cocycle condition and the ${\mathcal H}^{\alpha}_0$s define a
 {\em global} extension of ${\mathcal G}$ by ${\mathcal F}$.
  Essentially, this is the proof of the   fact that $Ext^1({\mathcal F},\,{\mathcal G})$  is
   the group of equivalence classes of length $1$ extensions of ${\mathcal G}$ by ${\mathcal F}$.

\smallskip

The proof of Theorem \ref{ThmH0Extp} can extended directly to the
vertex operators in $H^q(X,\,{\mathcal Ext}^p({\mathcal
F},\,{\mathcal G}))$. Thus, we have the following theorem.

\begin{Thm}\label{ThmHqExtp} Under the hypotheses of Theorem
\ref{ThmH0Extp}, each element of $H^q(X,\,{\mathcal
Ext}^p({\mathcal F},\,{\mathcal G}))$, with $p>0$ determines  a
$p$ extension of ${\mathcal G}_{|U_{\alpha_0\dots\alpha_q}}$ by
${\mathcal F}_{|U_{\alpha_0\dots\alpha_q}}$, where
$$U_{\alpha_0\dots\alpha_q}=\bigcap_{j=0}^qU_{\alpha_j}.$$
\end{Thm}

In the case of strings between two $D0$ branes the vertex
 operators define  global extensions; in fact, they coincide  with the respective
 %spaces of
  strings. In the following example we consider this point.

 \smallskip

\noindent
 {\it Example.} Let us assume that ${\mathcal I}\subset{\mathcal  O}$ is a sheaf of regular ideals such that
  $Z$, the support of ${\mathcal O}/{\mathcal I}$, has dimension
  $0$. We will study the space of strings stretching between ${\mathcal O}/{\mathcal
  I}$ and itself.  The set $Z$ can be regarded as a ringed space
  with structure sheaf ${\mathcal O}_Z={\mathcal O}/{\mathcal I}$.
  A free resolution of ${\mathcal O}_Z$ is the Koszul complex,
  that we recall briefly (see \cite{Eisenbud, G-H} for details).

  If $e_1,\dots,e_n$ is a set of symbols, we put ${\mathcal E}_p$
  for the free  ${\mathcal O}$-module generated by
  $\{e_{i_1}\wedge\dots\wedge e_{i_p}\}$, with $i_1< \dots
  <{i_p}$.
  That is, ${\mathcal E}_p$ is a free module with dimension $r:={n\choose p}$.

   The operator $\partial:{\mathcal E}_p\to {\mathcal
  E}_{p-1}$ is defined as follows: If $f_1,\dots, f_n$ is a regular sequence of functions which generate  the ideal
  ${\mathcal I}$ on an open $U$, then $\partial $ restricted to
  $U$ is defined by
  \begin{equation}\label{partial}
  \partial(e_{i_1}\wedge\dots\wedge
  e_{i_p})=\sum_{j=1}^p(-1)^{j-1}f_{i_j}e_{i_1}\wedge\dots\widehat{e}_{i_j}\dots\wedge
  e_{i_p}.
   \end{equation}
 The following exact sequence is a free resolution of ${\mathcal O}_Z$
 $$0\to{\mathcal E}_n\overset{\partial}{\to}{\mathcal E}_{n-1}
\overset{\partial}{\to}\dots\overset{\partial} {\to}{\mathcal
E}_0={\mathcal O}\overset{\rm proj.}{\longrightarrow} {\mathcal
 O}_Z\to 0.$$
 Applying the functor ${\mathcal Hom}(\,.\, ,\, {\mathcal O}_Z)$, one obtains
 the complex
 \begin{equation}\label{bigHom}
 \big({\mathcal Hom}({\mathcal E}_{\bullet},\,{\mathcal O}_Z),\,\partial ^*\big).
  \end{equation}

  As the sheaf ${\mathcal E}_p$ is the direct sum ${\mathcal O}^{\oplus\, r}$,
 \begin{equation}\label{mathcalHom}
 {\mathcal Hom}( {\mathcal E}_p,\, {\mathcal O}_Z)=\big({\mathcal O}_Z\big)^{\oplus\, r}.
 %{\oplus{n\choose p}}.
  \end{equation}
  Hence, the support of ${\mathcal Hom}( {\mathcal E}_p,\, {\mathcal O}_Z)$ is
  $Z$.
  On the other hand, since the $f_j$ vanish on $Z$, the  coboundary $\partial ^*$ operator
  induces the morphism zero in the stalk of
  ${\mathcal Hom}( {\mathcal E}_p,\, {\mathcal O}_Z)$  at any point of $Z$. Thus, the operator $\partial^*$ is identically
  zero.

The sheaves ${\mathcal Ext}^p({\mathcal O}_Z,\,{\mathcal O}_Z)$
are the homology elements of  the trivial complex (\ref{bigHom}).
So,
 \begin{equation}\label{D0D0}
{\mathcal Ext}^p({\mathcal O}_Z,\,{\mathcal O}_Z)= {\mathcal
Hom}( {\mathcal E}_p,\, {\mathcal O}_Z).
 \end{equation}

As ${\mathcal O}_Z$ is a skyscraper sheaf,
$$H^0(X,\,{\mathcal
O}_Z)=\bigoplus_{x\in Z}\big({\mathcal O}_Z\big)_x\simeq
\oplus_{x\in Z}{\mathbb C},$$
 and the other cohomology groups vanish. So, from (\ref{D0D0}) together with (\ref{mathcalHom}), it follows that
 \begin{align}H^q(X,\,{\mathcal Ext}^p({\mathcal O}_Z,\,{\mathcal O}_Z))&=0,\,\;\;\;\hbox{for}\, q\ne
 0. \notag \\
  H^0(X,\,{\mathcal Ext}^p({\mathcal O}_Z,\,{\mathcal
  O}_Z))&\simeq \oplus_{x\in Z}{\mathbb C}^r. \notag
 \end{align}

 Therefore, in the second page of the Local-to-Global spectral sequence all the rows are identically zero unless one.
 Thus,
 $$Ext^p({\mathcal O}_Z,\,{\mathcal O}_Z)=H^0(X,\,{\mathcal Ext}^p ( {\mathcal O}_Z,\,{\mathcal
 O}_Z)). $$
 That is, all the vertex operators for strings from the $D0$ brane
 ${\mathcal O}_Z$ to itself are {\em global} extensions.

\smallskip

 This  result, for the case $p=1$, can be also deduced from the observations explained in
 {\it Remark 1}.
  Let $\mathfrak{U}$ be a covering of $X$. Given $f\in H^0(X,\,{\mathcal Ext}^1({\mathcal O}_Z,\,{\mathcal O}_Z))$, it defines
 $f_{\alpha}\in \Gamma(U_{\alpha},\, {\mathcal Hom}({\mathcal E}_1,\,{\mathcal O}_Z))$, as in (\ref{cocyclefalpha}). Assuming that
    each point of the discrete set  $Z$ belongs to only one member of   ${\mathfrak
 U}$, then the spaces $\Gamma(U_{\alpha\beta},\,{\mathcal Hom}({\mathcal
 E}_0,\,{\mathcal O}_Z))$ and $\Gamma(U_{\alpha\beta},\,{\mathcal Hom}({\mathcal
 E}_1,\,{\mathcal O}_Z))$ are zero, since $U_{\alpha\beta}\cap
 Z=\emptyset$.  Hence taking $h=0\in\Gamma(U_{\alpha\beta},\,{\mathcal Hom}({\mathcal
 E}_0,\,{\mathcal O}_Z)),$  the conditions (\ref{globaxten}) are trivially satisfied, and $f$ determines an element of $Ext^1({\mathcal O}_Z,\,{\mathcal O}_Z)$.
 
%%%%%%%%%%%%%%%%%%%%%%%%%%%%%%%%%%%%%%%%%%%%%%%%%%%%%%%%%%%%%%%%%%%%%%%%%%%%%%%%%%%%%%%%%%%%%%%%%%%%%%%%%%%%%%%%%%%%%%%%%%%%%%%%%%%%%

%%%%%%%%%%%%%%%%%%%%%%%%%%%%%%%%%%%%%%%%%%%%%%%%%%%%%%%%%%%%%%%%%%%%%%%%%%%%%%%%%%%%%%%%%%%%%%%%%%%%%%%%%%%%%%%%%%%%%%%%%%%%%%%%%%%

\section{Final remarks.} We have applied methods and ideas from the cohomology of sheaves and from the theory of categories, for determining  properties of the spaces of strings between {\em general} branes, i.e., objects in the derived category of coherent sheaves. We have also defined the vertex operators for those strings. 

That generality had not been considered in the literature. The vertex operators usually considered are
 operators 
% In general, the study was limited to vertex operators 
for strings stretching between the sheaves defined by the sections of vector bundles.    Vertex operators for some Ext groups between non locally free sheaves are considered only exceptionally; for example, 
 % A exception is the paper  in \cite{K-Sharpe},
 by studying  branes on $X$ that are pushforward  of holomorphic vector bundles on a submanifold of $X$ \cite{K-Sharpe}. 
 The branes usually taken into account are sheaves defined by sections of vector bundles on $X$ or vector bundles defined on subspaces of $X$, as the ones that arise from the $K$-homology of $X$ \cite{Baum,Jia,R-S}. Thus, we present  applications of cohomological  methods to the study of spaces of vertex operators that had not yet been considered.

When $X$ is a manifold acted by a Lie group $G$, it seems natural to define ``equivariant" branes, and it is expected that the spaces of vertex operators for strings stretching between those branes support  representations of $G$. In the particular case, when $X$ is a flag manifold of the group $G$ and the branes are locally free equivariant sheaves, perhaps the Borel-Weil-Bott theorem allows us to characterize the corresponding space of vertex operators. We think that the   equivariance in the context of the  branes is an issue  worth exploring.  

In this article, we used mathematical procedures and results in order to better  understand concepts relative to the $B$-branes. In this context, the difficulty  to derive the results consists in translating the physical concepts to a mathematical language  and then the  application of  the appropriate mathematical tools. In general, these mathematical theories are far from the usual scientific background of  physicists.  

 One   ambitious different program  would be to consider possible suggestions and insights for Mathematics that can be derived from methods used by studying the physics of the branes.  This is an open tremendous    challenge (see \cite{MooreImpac}).

%The branes defined from the
 %$K$-homology of $X$ 

 % vertex operators for some Ext groups between non locally free sheaves In articles over $K$-homology

%Even though works with a strong mathematical deal with general branes
 
%%%%%%%%%%%%%%%%%%%%%%%%%%%%%%%%%%%%%
%That generality had not been considered in the literature; the branes usually considered are sheaves defined by sections of vector bundles on $X$ o vector bundles defined on subspaces of $X$, as the ones that arise from the $K$-homology of $X$ \cite{Jia},\cite{R-S}.

%%%%%%%%%%%%%%%%%%%%%%%%%%%%%%%%%%%%%%%%%%%%%%%%%%%%%%%%%%%%%%%%%%%%%%%%%%%%%%%%%%%%%%%%%%%%%%%%%%%%%%%%%%%%%%%%%%%%%%%%%%
%%%%%%%%%%%%%%%%%%%%%%%%%%%%%%%%%%%%%%%%%%%%%%%%%%%%%%%%%%%%%%%%%%%%%%%%%%%%%%%%%%%%%%%%%%%%%%%%%%%%%%%%%%%%%%%%%%%%%
\section{Appendix}

In this section, we will prove Proposition \ref{P:Appendix}.

 Let $K,L$ be complexes in an additive category $\mathfrak{A}$ and $f$ a morphism from $K$ to $L$. We denote
  by ${\rm Con}(f)$ and ${\rm Cyl}(f)$ the corresponding mapping cone and mapping cylinder complexes.
 Fixed a complex $F$, we put $\widehat L$ for denoting the complex of
abelian groups
  ${\rm Hom}^{\bullet}_\mathfrak{A}(F,L)$,   analogously $\widehat K={\rm Hom}^{\bullet}_\mathfrak{A}(F,K)$ and
  $\widehat f:= {\rm Hom}^{\bullet}_\mathfrak{A}(F,f)$.
  %In the Appendix we will give a proof the following result.
 \begin{Lem}\label{L:Cone-Cylinder}
 $${\rm Con}(\widehat f)=\widehat{{\rm Con}(f)},\;\; {\rm Cyl}(\widehat f)=\widehat{{\rm Cyl}(f)}.$$
 \end{Lem}

{\it Proof.} By definition
$$\widehat L^m=\prod_i{\rm Hom}_\mathfrak{A}(F^i,\,L^{m+i}),$$
and coboundary operator
\begin{equation}\label{dHatL}
d^m_{\widehat L}(s_i)=\big(d^{m+i}_{L}\circ s_i -(-1)^{m}s_{i+1}\circ
d^î_F \big),
\end{equation}
with $s_i\in{\rm Hom}_\mathfrak{A}(F^i,\,L^{m+i})$. The complex
$\widehat K$ and the operator $d_{\widehat K}$ are defined analogously.

As $f:K\to L$ is a morphism of complexes, it induces a complex
morphism $\widehat f:\widehat K\to \widehat L$ by the formula
\begin{equation}\label{Hatf2}
\widehat f(r_i)=(f^{m+1+i}\circ r_i),
\end{equation}
with $r_i\in{\rm Hom}_\mathfrak{A}(F^i,\,K^{m+1+i}).$

By definition,
\begin{equation}\label{ConeHatf}
{\rm Con}(\widehat f)^m=\widehat K^{m+1}\oplus \widehat L^m=\prod_i\Big( {\rm
Hom}_\mathfrak{A}(F^i,\,K^{m+1+i}) \oplus {\rm
Hom}_\mathfrak{A}(F^i,\,L^{m+i}) \Big).
 \end{equation}
 The coboundary operator of this complex is defined  by (see \cite[page
 23]{Iversen})
 \begin{equation}\label{cobound-Cone}
 d_{{\rm Con}(\widehat f)}=\begin{pmatrix}
-d_{\widehat K} & 0 \\
-\widehat f & d_{\widehat L}
 \end{pmatrix}.
  \end{equation}
That is,
 \begin{equation}\label{alignCone}
 d^m_{{\rm Con}(\widehat f)}\begin{pmatrix}
  r_i  \\
 s_i
  \end{pmatrix} =\begin{pmatrix} -d_K^{m+1+i}\circ r_i-(-1)^mr_{i+1}\circ d^i_F \\
  
-f^{m+1+i}\circ r_i+ d_L^{m+i}\circ s_i -(-1)^{m}s_{i+1}\circ d^i_F 
 \end{pmatrix}
 \end{equation}

%\begin{equation}\label{alignCone}
%\begin{align}\label{alignCone}
%&d^m_{{\rm Con}(\widehat f)}(r_i,\,s_i)=\\
%& \big(-d_K^{m+1+i}\circ r_i-(-1)^mr_{i+1}\circ d^i_F,\,
%-f^{m+1+i}\circ r_i+ d_L^{m+i}\circ s_i -(-1)^{m}s_{i+1}\circ d^i_F \big) \notag.
%\end{align}
%\end{equation}

 On the other hand, we consider the complex
 $$\widehat{{\rm Con}(f)}={\rm Hom}^{\bullet}_{\mathfrak A}(F,\,{\rm Con}(f)).$$
 Hence,
 $$\widehat{{\rm Con}(f)}^m=\prod_i{\rm Hom}_{\mathfrak A}(F^i,\, {\rm Con}(f)^{m+i})=\prod_i{\rm Hom}_{\mathfrak A}(F^i,\, K^{m+1+i}\oplus L^{m+i}).$$
 From (\ref{ConeHatf}), it follows $\widehat{{\rm Con}(f)}^m={\rm Con}(\widehat f)^m.$

 According to (\ref{dHatL}), the coboundary operator of $\widehat{{\rm Con}(f)}$ is defined by
 $$d^m_{\widehat{{\rm Con}(f)}}(r_i,\,s_i)=d^{m+i}_{{\rm Con}(f)}\circ(r_i,\,s_i)-(-1)^m(r_{i+1},\,s_{i+1})\circ d^i_F.$$
 That is,
 % $$d^m_{\widehat{{\rm Con}(f)}}(r_i,\,s_i)=
 %  \begin{pmatrix}
%-d_{ K} & 0 \\
%- f & d_L
% \end{pmatrix}
% \begin{pmatrix}
%  r_i  \\
%  s_i
%  \end{pmatrix}-(-1)^m(r_{i+1},\,s_{i+1})\circ d^i_F.$$

 $$d^m_{\widehat{{\rm Con}(f)}}\begin{pmatrix}
  r_i  \\
 s_i
  \end{pmatrix}=
  \begin{pmatrix}
-d_{ K} & 0 \\
- f & d_L
 \end{pmatrix}
 \begin{pmatrix}
  r_i  \\
 s_i
  \end{pmatrix}-(-1)^m\begin{pmatrix}
  r_{i+1}  \\
 s_{i+1}
  \end{pmatrix}\circ d^i_F.$$

  This expression coincides with (\ref{alignCone}).
  So, we have the equality of the complexes
 % $$\big( \widehat{{\rm Con}(f)},\,d_{\widehat{{\rm Con}(f)}}   \big)= \big({\rm Con}(\widehat f)},\, d_{{\rm Con}(\widehat f)}  \big).$$
$$\big(  \widehat{{\rm Con}(f)}  ,\, d_{\widehat{{\rm Con}(f)}}
\big)= \big( {\rm Con}(\widehat f),\, d_{{\rm Con}(\widehat f)}  \big).$$

  The proof of the property of relative to the cylinder is similar.

   \qed

In general, if $l:{\mathcal M}\to {\mathcal N}$ is a morphism in
the category ${\mathfrak C}$ of complexes of ${\mathcal
O}$-modules and $U$ an open subset of $X$, then
 %Let ${\mathcal A}$, ${\mathcal
%F}$, ${\mathcal G}$  be objects  in the category $\mathfrak{C}$ of
%complexes
 %of  ${\mathcal O}$-modules and    $h:{\mathcal A}\to{\mathcal G}$  a morphism  in this category. If
  % $U$ an open subset of $X$, then
\begin{equation}\label{Coneh-U}
{\rm Con}(l_{|\,U})={\mathcal M}_{|\,U}[1]\oplus{\mathcal
N}_{|\,U} = \big({\mathcal M}[1]\oplus{\mathcal N}\big)_{|\,U}
={\rm Con}(l)_{|\,U}.
 \end{equation}
  A  similar relation holds for the mapping cylinder.

\smallskip

   Let ${\mathcal P} $, ${\mathcal B} $, ${\mathcal A} $ objects of the category ${\mathfrak
   C}$,
   $g:{\mathcal B} \to{\mathcal A}$ a morphism in this category,
   and
   $U$ an open subset of $X$. Then $g_{|\,U}:{\mathcal B}_{|\,U}\to{\mathcal
   A}_{|\,U}$ is a morphism in the category of complexes of ${\mathcal
   O}_{|\,U}$-modules.

We put
$$\widehat g={\mathcal Hom}^{\bullet}({\mathcal P},\,g):{\mathcal Hom}^{\bullet}({\mathcal P},\,{\mathcal B})\to {\mathcal Hom}^{\bullet}({\mathcal P},\,{\mathcal
A}).$$
 Analogously,
  \begin{equation}\label{Hatg-U}
  \widehat{g_{|\,U}}:={Hom}^{\bullet}_{{\mathcal O}_{|\,U}}({\mathcal P}_{|\,U},\,g_{|\,U}):
  {Hom}^{\bullet}_{{\mathcal O}_{|\,U}}({\mathcal P}_{|\,U},\,{\mathcal B}_{|\,U})\to
  {Hom}^{\bullet}_{{\mathcal O}_{|\,U}}({\mathcal P}_{|\,U},\,{\mathcal
A}_{|\,U}).
 \end{equation}

\begin{Lem}\label{L:ConeHatU}
% $${\rm Con}(\widehat{g}_{|\,U})={\rm Con}(\widehat{g_{|\,U}}),\;\;{\rm Cyl}(\widehat{g}_{|\,U})={\rm Cyl}(\widehat{g_{|\,U}})$$
 $${\rm Con}(\widehat{g})(U)={\rm Con}(\widehat{g_{|\,U}}),\;\;
 {\rm Cyl}(\widehat{g})(U)= {\rm Cyl}(\widehat{g_{|\,U}}).$$
 \end{Lem}

 {\it Proof.} By definition of the cone mapping and the ${\mathcal Hom}$ functor,
\begin{align}
{\rm Con}(\widehat{g})(U) & = \big({\mathcal Hom}^{\bullet}({\mathcal P},\,{\mathcal B})[1]
\oplus{\mathcal Hom}^{\bullet}({\mathcal P},\,{\mathcal
A})\big)(U) \notag \\
 &={Hom}^{\bullet}_{{\mathcal O}_{|\,U}}({\mathcal P}_{|\,U},\,{\mathcal
 B}_{|\,U})[1]\oplus
{Hom}^{\bullet}_{{\mathcal O}_{|\,U}}({\mathcal
P}_{|\,U},\,{\mathcal A}_{|\,U}) \notag.
 \end{align}
 According (\ref{Hatg-U}), the last expression is precisely ${\rm Con}(\widehat{g_{|\, U}})$.

%On the other hand,
%\begin{equation} \label{Hatg_}
%\widehat g_{|\,U}:{\mathcal Hom}^{\bullet}({\mathcal P},\,{\mathcal
%B})_{|\,U}\to {\mathcal Hom}^{\bullet}({\mathcal P},\,{\mathcal
%A})_{|\,U}.
 %\end{equation}
% Thus,
% \begin{equation}\label{ConHatg_}
% {\rm Con}(\widehat g_{|\,U})={\mathcal Hom}^{\bullet}({\mathcal P},\,{\mathcal
% B})_{|\,U}[1]\oplus
%{\mathcal Hom}^{\bullet}({\mathcal P},\,{\mathcal A})_{|\,U}.
% \end{equation}
%From (\ref{Hatg_}) and (\ref{ConHatg_}), it follows ${\rm
%Con}(\widehat{g}_{|\,U})={\rm Con}(\widehat{g_{|\,U}}).$
 \qed

\smallskip

\noindent {\bf Proof of Proposition \ref{P:Appendix}.} By
definition of
 %the functor
 ${\mathcal Hom}$
\begin{equation} \label{P:AppendixProof1}
{\mathcal Hom}^{\bullet}\big( {\mathcal P},\,{\rm Con}(g)
\big)(U)= Hom^{\bullet}_{{\mathcal O}_{|\,U}} \big( {\mathcal
P}_{|\,U},\,{\rm Con}(g)_{|\,U}\big).
 \end{equation}
 By (\ref{Coneh-U})
\begin{equation}\label{anterior}
Hom^{\bullet}_{{\mathcal O}_{|\,U}}\big( {\mathcal
P}_{|\,U},\,{\rm Con}(g)_{|\,U}\big)=
 Hom^{\bullet}_{{\mathcal O}_{|\,U}}
\big( {\mathcal P}_{|\,U},\,{\rm Con}(g_{|\,U})\big).
 \end{equation}

  Denoting
${\mathfrak A}$ the category of complex of ${\mathcal
O}_{|\,U}$-modules, then $Hom^{\bullet}_{{\mathcal O}_{|\,U}}={\rm
Hom}^{\bullet}_{\mathfrak A}.$ With the notation used at the
beginning  of the Appendix,
%By Lemma \ref{L:Cone-Cylinder},
the right hand side of (\ref{anterior}) is $\widehat{{\rm
Con}(g_{|\,U}}).$
 From Lemmas \ref{L:Cone-Cylinder} and
 \ref{L:ConeHatU},
 %together with (\ref{Coneh-U}),
  it follows
\begin{equation}\label{HatCongU}
\widehat{{\rm Con}(g_{|\,U}})={\rm Con}(\widehat{g_{|\,U}})
%{\rm Con}(\widehat{g}_{|\,U})
={\rm Con}(\widehat g )(U).
 \end{equation}
 From (\ref{P:AppendixProof1}), together with (\ref{anterior}) and (\ref{HatCongU}), we
 obtain
$${\mathcal Hom}^{\bullet}({\mathcal P},\,{\rm Con}(g))(U)={\rm
Con}(\widehat g)(U).$$

 Similarly,
 \begin{equation}\label{CylinderHat}
 {\mathcal Hom}^{\bullet}({\mathcal P},\,{\rm Cyl}(g))(U)={\rm Cyl}(\widehat g)(U).
 \end{equation}

\qed

%%%%%%%%%%%%%%%%%%%%%%%%%%%%%%%%%%%%%%%%%%%%%%%%%%%%%%%%%%%%%%%%%%%%%%%%%%%%%%%%%%%%%%%%%%%%%%%%%%%%%%%%%%%%%%%%%%%%%%%%%%%%%%%%%%%%%%
%%%%%%%%%%%%%%%%%%%%%%%%%%%%%%%%%%%%%%%%%%%%%%%%%%%%%%%%%%%%%%%%%%%%%%%%%%%%%%%%%%%%%%%%%%%%%%%%%%%%%%%%%%%%%%%%%%%%%%%%%%%%%%%%%%%%%%%%%


\begin{thebibliography}{99}





 \bibitem{Aspin}
  Aspinwall, P. S.: $D$-branes on Calabi-Yau manifolds.
  In Progress in String Theory. Pages 1-152. World Sci. Publ.
  %Hackensack, NJ.
  (2005).




\bibitem{Aspin-et}
Aspinwall, P. S. et al.: Dirichlet branes and mirror symmetry. Clay
mathematics monographs vol 4. Amer. Math. Soc. (2009).



%\bibitem{Aspin-Law}
%Aspinwall, P. S., Lawrence, A. E.: Derived categories and
%zero-brane stability. JHEP {\bf 08}, 004 (2001).

\bibitem{Baum}
 Baum, P.: $K$-homology and $D$-branes. In Superstrings, Geometry, Topology, and $C^*$-algebras. Vol 81, Proc. Sympos. Pure Math. 81-94. AMS (2009).


\bibitem{B-K}
Bondal A., Kapranov M.: Representable functors, Serre functors,
and mutations, Izv. Akad. Nauk SSSR, Ser.Mat., {\bf 53}, 1183-1205
(1989); English transl. in Math. USSR Izv., {\bf 35}, 519-541
(1990).




%\bibitem{Borelli}
%Borelli, M.: Some results on ampleness and divisorial  schemes.
%Pacific J. Math. {\bf 23}, 217-227 (1967).



%\bibitem{Douglas}
% Douglas, M. R.: D-branes, categories and $N=1$ supersymmetry. J. Math. Phys {\bf 42}, 2818-2843 (2001).

\bibitem{Fulton}
 Fulton, W.: Intersection theory. Springer
 (1998).

\bibitem{Eisenbud}
Eisenbud, E.: Commutative algebra with a view toward algebraic
geometry. Springer-Verlag (1995).

\bibitem{Ge-Ma}
Gelfand, S. I., Manin, Y. I.: Methods of homological algebra. Springer (2003).

\bibitem{G-H}
 Griffiths, P.,  Harris, J.: Principles of algebraic geometry.
 John Wiley (1994).

\bibitem{Hart}
  Hartshorne,  R.:
   Algebraic  geometry. Springer-Verlag (1983).

%\bibitem{Harvey}
% Harvey, J. A.: TASI 2003 Lectures on anomalies. 2005.


\bibitem{Iversen}
Iversen, B.: Cohomology of sheaves.  Springer (1986).

\bibitem{Jia}
 Jia, B.: D-branes and $K$-homology. arXiv:1306.0535\,[math. KT].



\bibitem{K-Sharpe}
Katz, S., Sharpe, E.: D-branes, open string vertex operators, and Ext groups.
 Adv. Theor. Math. Phys. {\bf 6}, 979-1030 (2003).

\bibitem{Kas-Sch}
 Kashiwara, M., Schapira, P.:  Sheaves on manifolds.
  Springer-Verlag (2002).




%\bibitem{Koszul}
 %Koszul, J. L.: Formes hermitiennes canoniques des espaces
% homog\`enes complexes.
% S\'eminaire N. Bourbaki, 1954-1956, exp. 108, 69-75.

%\bibitem{Lam}
%Lam, T. Y.: Lectures on modules and rings. Springer-Verlag (1999).



\bibitem{M-L}
Mac Lane, S.: Homology. Springer-Verlag (1975).


%\bibitem{M-M}
% Minasian, R.,  Moore, G.: K-theory and Ramond-Ramond charge. J High
% Energy Phys. {\bf 11}, (1997) 002.

\bibitem{Mitchell}
 Mitchell, B.: Theory of categories. Academic Press. (1965).

\bibitem{MooreImpac}
 Moore, G. W.:          The impact of $D$-branes in Mathematics. (2014) (In www.physics.rutgers.edu/$\sim$gmoore).
%\bibitem{Oda}
% Oda, T.: Convex bodies and algebraic geometry. Springer-Verlag (1988).

\bibitem{R-S}
Reis, R.M.G.,  Szabo, R.J.: Geometric K-homology of flat D-branes.  Commun. Math. Phys. {\bf 266}, 71-122 (2006).

%\bibitem{R}
%Rotman, J. J.: An introduction to homological algebra. Springer
%(2008).

\bibitem{Serre}
Serre, J. P.: Faisceaux alg\'ebriques coh\'erents. Annals of Math.
{\bf 61}, 197-278 (1955).






%\bibitem{V4}
 %  Vi\~na,  A.:  Equivariant branes.


\bibitem{Wei}
Weibel, Ch. A.: An introduction to homological algebra. Cambridge
U.P. (1997).




\bibitem{Witten0}
Witten, E.: Chern-Simmons gauge theory as a string theory. In H.Hofer et alt. editors The Floer memorial volume. 637-678. Birkh\"auser (1995).

\bibitem{Witten}
Witten, E.: Mirror manifolds and topological field theory.
 In Mirror Symmetry I. Edit S.-T. Yau.  Pages 121-157. AMS (1998).





\end{thebibliography}
\end{document}